\documentclass[notitlepage]{report}
\usepackage{amsmath}
\usepackage{amstext}
\usepackage{amsthm}
\usepackage{amsxtra}
\usepackage{amssymb}
\usepackage{amsfonts}
\usepackage{epsfig}
\usepackage{float}
\usepackage{cite}
\usepackage{rotating}
\usepackage{color}
\usepackage{lscape}
\usepackage{subfigure}
\usepackage{graphicx}
\usepackage{empheq}
\usepackage{caption2}
\usepackage[section]{placeins}
\usepackage{indentfirst}
\usepackage{rrprefs}
\title{\bf The polynomial method for random matrices}
\author{\sc N. Raj Rao\thanks{%
        MIT Department of Electrical Engineering and Computer Science,
        Cambridge, MA 02139,
        raj@mit.edu}
        \and 
        \sc Alan Edelman\thanks{%
        MIT Department of Mathematics,
        Cambridge, MA 02139,
        edelman@math.mit.edu}}
\date{\small\sc \today}
\begin{document}
\maketitle
\captionstyle{hang}
\begin{abstract}

We define a class of ``algebraic'' random matrices. These are random matrices for which the Stieltjes transform of the limiting eigenvalue distribution function is algebraic, \ie, it satisfies a (bivariate) polynomial equation. The Wigner and Wishart matrices whose limiting eigenvalue distributions are given by the semi-circle law and the Mar\v{c}enko-Pastur law are special cases. 
  
Algebraicity of a random matrix sequence is shown to act as a certificate of the computability of the limiting eigenvalue density function. The limiting moments of algebraic random matrix sequences, when they exist, are shown to satisfy a finite depth linear recursion so that they may often be efficiently enumerated in closed form.

In this article, we develop the mathematics of the {\em polynomial method} which allows us to describe the class of algebraic matrices by its generators and map the constructive approach we employ when proving algebraicity into a software implementation that is available for download in the form of the RMTool random matrix ``calculator'' package. Our characterization of the closure  of algebraic probability distributions under free additive and multiplicative convolution operations allows us to simultaneously establish a framework for computational (non-commutative) ``free probability'' theory. We hope that the tools developed allow researchers to finally harness the power of the infinite random matrix theory.
\end{abstract}

\pagestyle{myheadings}
\thispagestyle{plain}
\markboth{N. R. Rao and  A. Edelman}{The polynomial method}

\section{Introduction}
\begin{figure}
\centering
\subfigure[The limiting eigenvalue density function for the GOE and Wishart matrices.]{
\label{fig:wishart example}
\includegraphics[width=3.75in]{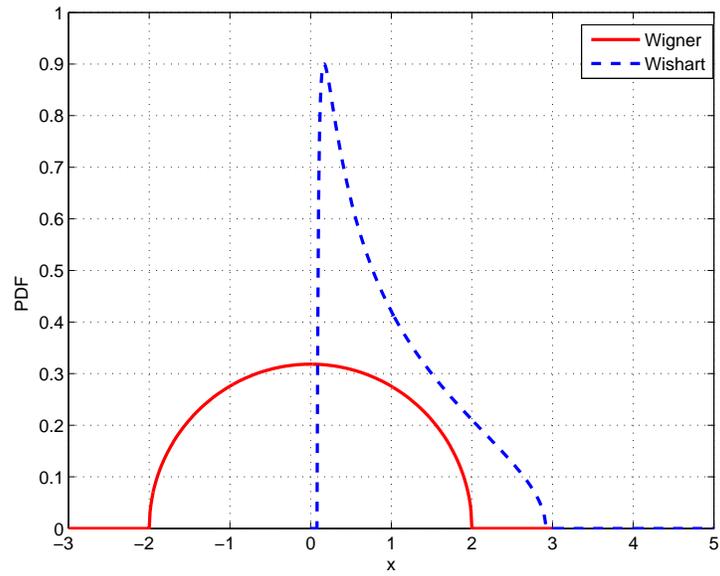}
}\\
%\vspace{-0.2in}
\subfigure[The limiting eigenvalue density function for the sum and product of independent GOE and Wishart matrices.]{
\label{fig:wigpluswis example}
\includegraphics[width=3.75in]{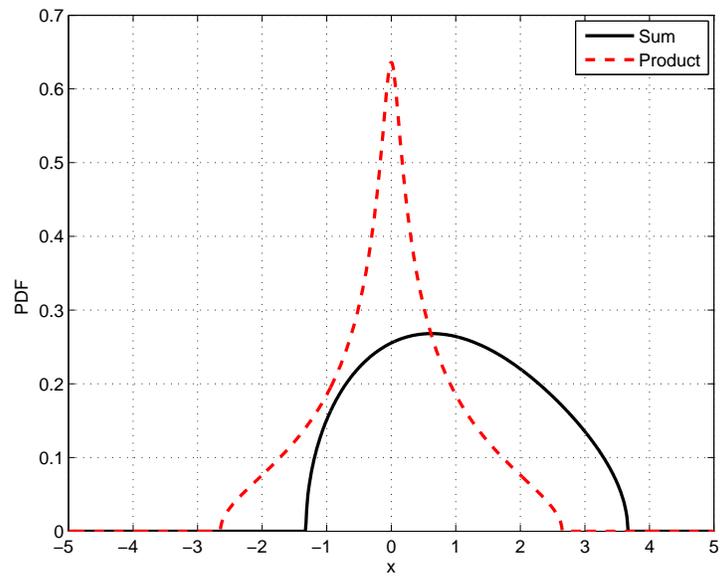}
}\\
\vspace{-0.1in}
\caption{A representative computation using the random matrix calculator.}
\label{fig:calculator}
\end{figure}

We propose a powerful method that allows us to calculate the limiting eigenvalue distribution  of a large class of random matrices. We see this method as allowing us to expand our reach beyond the well known special random matrices whose limiting eigenvalue distributions have the semi-circle density \cite{wigner55a}, the Mar\v{c}enko-Pastur density \cite{marcenko67a}, the McKay density \cite{mckay81a} or their close cousins \cite{silverstein85a,collins05a}. 
In particular, we encode transforms of the limiting eigenvalue distribution function as solutions of bivariate polynomial equations.  Then canonical operations on the random matrices become operations on the bivariate polynomials. 
We illustrate this with a simple example. Suppose we take the Wigner matrix, sampled in \matlab as: 

\begin{small}
\begin{verbatim}
     G = sign(randn(N))/sqrt(N); A = (G+G')/sqrt(2);
\end{verbatim}
\end{small}
\noindent whose eigenvalues in the $N \to \infty$ limit follow the semicircle law, and the Wishart matrix which may be sampled in \matlab as:
\begin{small}
\begin{verbatim}
     G = randn(N,2*N)/sqrt(2*N); B = G*G';
\end{verbatim}
\end{small}

\noindent whose eigenvalues in the limit follow the Mar\v{c}enko-Pastur law. The associated limiting eigenvalue distribution functions have Stieltjes transforms $m_{A}(z)$ and $m_{B}(z)$ that are solutions of the equations $\lmz{A}=0$ and $\lmz{B}=0$, respectively,  where $$\lmz{A}= m^2+z\,m+1, \qquad \lmz{B} = {m}^{2}z- \left( -2\,z+1 \right) m+2.$$ 
The sum and product of independent samples of these random matrices have limiting eigenvalue distribution functions whose Stieltjes transform is a solution of the bivariate polynomial equations $\lmz{A+B}=0$ and $\lmz{AB}=0$, respectively, which can be calculated from $\lmzs{A}$ and $\lmzs{B}$ alone.
To obtain $\lmz{A+B}$ we apply the transformation labelled as ``Add Atomic Wishart'' in Table \ref{tab:summary bivariate transforms} with $c = 2$, $p_{1} = 1$ and $\lambda_{1} = 1/c = 0.5$ to obtain the operational law
\begin{equation}\label{eq:intro ex add}
\lmz{A+B} = \lmzs{A}\left(m,z-\frac{1}{1+0.5m}\right). 
\end{equation}
Substituting  $\lmzs{A}=m^2+z\,m+1$ in (\ref{eq:intro ex add}) and clearing the denominator, yields the bivariate polynomial
\begin{equation}
\lmz{A+B}={m}^{3}+ \left( z+2 \right) {m}^{2}- \left( -2\,z+1 \right) m+2.
\end{equation} 
Similarly, to obtain $\lmzs{AB}$, we apply the transformation labelled as ``Multiply Wishart'' in Table \ref{tab:summary bivariate transforms}
with $c = 0.5$ to obtain the operational law
\begin{equation}\label{eq:intro ex prod}
\lmz{AB} = \lmzs{A}\left((0.5-0.5zm)\,m,\dfrac{z}{0.5-0.5zm}\right). 
\end{equation}
Substituting $\lmzs{A}=m^2+z\,m+1$ in (\ref{eq:intro ex prod}) and clearing the denominator, yields the bivariate polynomial
\begin{equation}
\lmz{AB}= {m}^{4}{z}^{2}-2\,{m}^{3}z+{m}^{2}+4\,mz+4.
\end{equation}
Figure \ref{fig:calculator} plots the density function associated with the limiting eigenvalue distribution for the Wigner and Wishart matrices as well as their sum and product extracted directly from $\lmz{A+B}$ and $\lmz{AB}$. In these examples, algebraically extracting the roots of these polynomials using the cubic or quartic formulas is of little use except to determine the limiting density function. As we shall demonstrate in Section \ref{sec:moments}, the algebraicity of the limiting distribution (in the sense made precise next) is what allows us to readily enumerate the moments efficiently directly from the polynomials $\lmz{A+B}$ and $\lmz{AB}$.
\subsection{Algebraic random matrices: Definition and Utility}
A central object in the study of large random matrices is the empirical distribution function which is defined, for an $N \times N$ matrix ${\bf A}_{N}$ with real eigenvalues, as 
\begin{equation}\label{eq:intro edf}
F^{{\bf A}_{N}}(x)=\frac{\textrm{Number of eigenvalues of }{\bf A}_{N} \leq x }{N}.
\end{equation}
For a large class of random matrices, the empirical distribution function $F^{{\bf A}_{N}}(x)$ converges, for every $x$, almost surely (or in probability) as $N \to \infty$ to a non-random distribution function $F^{A}(x)$. The dominant theme of this paper is that ``algebraic'' random matrices form an important subclass of analytically tractable random matrices and can be effectively studied using combinatorial and analytical techniques that we bring into sharper focus in this paper.

\theorembox{
\begin{definitionS}[Algebraic random matrices]\label{def:intro Malg}
Let $F^{A}(x)$ denote the limiting eigenvalue distribution function of a sequence of random matrices ${\bf A}_{N}$. If a bivariate polynomial $L_{{\rm mz}}(m,z)$ exists such that
\[
m_{A}(z) = \int \frac{1}{x-z} dF^{A}(x) \qquad z \in \mathbb{C}^{+} \setminus \mathbb{R}
\]
is a solution of $L_{{\rm mz}}(m_{A}(z),z) = 0$ then ${\bf A}_{N}$ is said to be an algebraic random matrix. The density function $f_{A} :=dF^{A}$ (in the distributional sense) is referred to as an algebraic density and we say that ${\bf A}_{N} \in \mathcal{M}_{{\rm alg}}$, the class of algebraic random matrices and $f_{A} \in \Palg$, the class of algebraic distributions.
\end{definitionS}
}

\begin{figure}[t]
%%%\captionstyle{flushleft}
\centering
\begin{center}
\setlength{\unitlength}{2ex}
\begin{picture}(49,29)
%%% title
%%\put(12.25,-3){\Large\bf A Random Matrix Calculator}
%%% left graph (A)
\put(16,22){$A$}
\put(5,14){\vector(0,1){9}}
\put(5,14){\vector(1,0){12}}
\qbezier(7,14)(7,22)(8,22)
\qbezier(8,22)(15,18)(15,14)
\put(10.5,15){\makebox(0,0)[b]{\begin{tabular}{l}
    Limiting\\ Density\end{tabular}}}
%%% box around (A)
\multiput(4,13)(14,0){2}{\line(0,1){11}}
\multiput(4,13)(0,11){2}{\line(1,0){14}}
%%% arrow between (A) and (B)
\put(19,18){\vector(1,0){2}}
%%% right graph (B)
\put(34,22){$B$}
\put(23,14){\vector(0,1){9}}
\put(23,14){\vector(1,0){12}}
\qbezier(25,14)(25,22)(26.5,22)
\qbezier(26.5,22)(27,22)(27.5,21)
\qbezier(27.5,21)(28,20)(29,20)
\qbezier(29,20)(30,20)(30.5,20.5)
\qbezier(30.5,20.5)(31,21)(31.5,21)
\qbezier(31.5,21)(33,21)(33,14)
\put(29,15){\makebox(0,0)[b]{\begin{tabular}{c}
    Limiting\\ Density\end{tabular}}}
%%% box around (B)
\multiput(22,13)(14,0){2}{\line(0,1){11}}
\multiput(22,13)(0,11){2}{\line(1,0){14}}
%%% Deterministic
\put(4.5,8.5){$A+\alpha\,I$}
\put(14,8.5){$\alpha\times A$}
\put(23.5,8.5){$A^{-1}$}
\put(31.3,8.5){$\displaystyle\frac{p A+q I}{r A+s I}$}
\multiput(3,7)(9,0){4}{\line(0,1){4}}
\multiput(10,7)(9,0){4}{\line(0,1){4}}
\multiput(3,11)(9,0){4}{\line(1,0){7}}
\multiput(3,7)(9,0){4}{\line(1,0){7}}
%\put(40.5,9){\vector(-1,0){2}}
%\put(41,8.5){Deterministic}
\put(20,12){\makebox(0,0){Deterministic}}
%%% Stochastic
\put(4,2.5){$A+W(c)$}
\put(13,2.5){$W(c)\times A$}
\put(21.5,2.5){$W^{-1}(c) \times A$}
\put(30.25,2.6){\small{$\begin{array}{c} (A^{1/2}+G)\\[-5pt] \times\\[-0.75pt]
    (A^{1/2}+G)'\end{array}$}}
\multiput(3,1)(9,0){4}{\line(0,1){4}}
\multiput(10,1)(9,0){4}{\line(0,1){4}}
\multiput(3,5)(9,0){4}{\line(1,0){7}}
\multiput(3,1)(9,0){4}{\line(1,0){7}}
%\put(40.5,3){\vector(-1,0){2}}
%\put(41,2.5){Stochastic}
\put(20,6){\makebox(0,0){Stochastic}}
%%% box it all in
\multiput(2,0)(36,0){2}{\line(0,1){25}}
\multiput(2,0)(0,25){2}{\line(1,0){36}}
%%%\put(14.5,0){\Large\bf A Random Matrix Calculator}
\end{picture}
\end{center}
\caption{A random matrix calculator where a sequence of deterministic and stochastic operations performed on an algebraic random matrix sequence ${\bf A}_{N}$ produces an algebraic random matrix sequence ${\bf B}_{N}$. The limiting eigenvalue density and moments of a algebraic random matrix can be computed numerically, with the latter often in closed form.}
\label{fig:rmtool}
\end{figure}
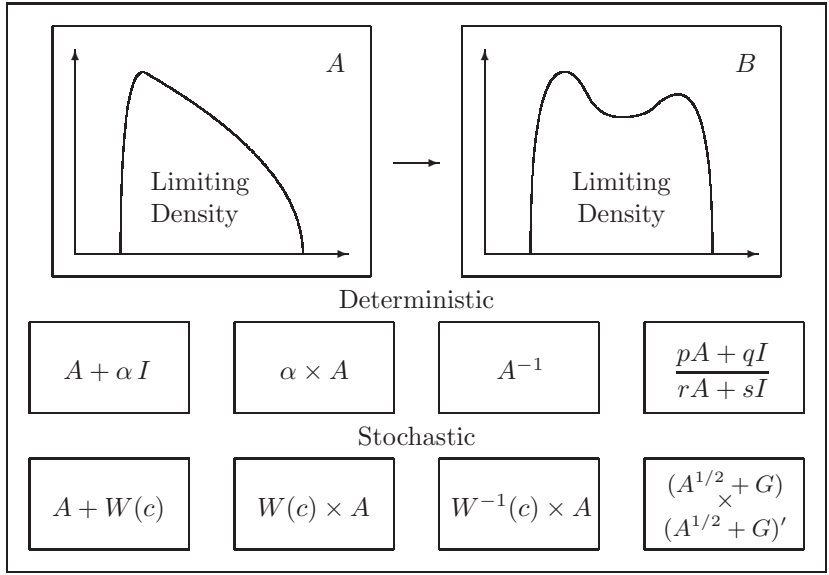

The utility of this, admittedly technical, definition comes from the fact that we are able to concretely specify the generators of this class. We illustrate this with a simple example. Let ${\bf G}$ be an $n \times m$ random matrix with i.i.d. standard normal entries with variance $1/m$. The matrix ${\bf W}(c)  = {\bf G}{\bf G}'$ is the Wishart matrix parameterized by $c = n/m$. Let ${\bf A}$ be an arbitrary algebraic random matrix independent of ${\bf W}(c)$.  
Figure \ref{fig:rmtool} identifies deterministic and stochastic operations that can be performed on ${\bf A}$ so that the resulting matrix is algebraic as well. The calculator analogy is apt because once we start with an algebraic random matrix, if we keep pushing away at the buttons we still get an algebraic random matrix whose limiting eigenvalue distribution is concretely computable using the algorithms developed in Section \ref{sec:pol method computational aspects}. 

The algebraicity definition is important because everything we want to know about the limiting eigenvalue distribution of ${\bf A}$ is encoded in the bivariate polynomial $\lmz{A}$. In this paper, we establish the algebraicity of each of the transformations in Figure \ref{fig:rmtool} using the ``hard'' approach that we label as the {\em polynomial method} whereby we explicitly determine the operational law for the polynomial transformation $\lmz{A} \mapsto \lmz{B}$ corresponding to the random matrix transformation ${\bf A} \mapsto {\bf B}$.  This is in contrast to the ``soft'' approach taken in a recent paper by Anderson and Zeitouni \cite[Section 6]{anderson06c} where the algebraicity of Stieltjes transforms under hypotheses frequently fulfilled in RMT is proven using dimension theory for noetherian local rings. 
The catalogue of admissible transformations, the corresponding ``hard'' operational law and their software realization is found in Section \ref{sec:pol method computational aspects}. This then allows us to calculate the eigenvalue distribution functions of a large class of algebraic random matrices that are generated from other algebraic random matrices. 
In the simple case involving Wigner and Wishart matrices considered earlier, the transformed polynomials were obtained by hand calculation. Along with the theory of algebraic random matrices we also develop a software realization that maps the entire catalog of transformations (see Tables \ref{tab:summary bivariate transforms} -\ref{tab:remaining transformations}) into symbolic \matlab code. Thus, for the same example, the sequence of commands:

\begin{small}
\begin{verbatim}
>> syms m z
>> LmzA = m^2+z*m+1;
>> LmzB = m^2-(-2*z+1)*m+2;
>> LmzApB = AplusB(LmzA,LmzB);
>> LmzAtB = AtimesB(LmzA,LmzB);
\end{verbatim}
\end{small}

\noindent could also have been used to obtain $\lmzs{A+B}$ and $\lmzs{AB}$. We note that the commands \texttt{AplusB} and \texttt{AtimesB} implicitly use the free convolution machinery (see Section \ref{sec:computational free probability}) to perform the said computation.
To summarize, by defining the class of algebraic random matrices, we are  able to extend the reach of infinite random matrix theory well beyond the special cases of matrices with Gaussian entries. The key idea is that by encoding probability densities as solutions of bivariate polynomial equations, and deriving the correct operational laws on this encoding, we can take advantage of powerful symbolic and numerical techniques to compute these densities and their associated moments. 
\subsection{Outline}
This paper is organized as follows. We introduce various transform representations of the distribution function in Section \ref{sec:transforms}.  We define algebraic distributions and the various manners in which they can be implicitly  represented in \ref{sec:alg distributions} and describe how they may be algebraically  manipulated in \ref{sec:kronecker}. 
The class of algebraic random matrices is described in Section \ref{sec:algebraic random matrices} where the theorems are stated and proved by obtaining the operational law on the bivariate polynomials summarized in Section \ref{sec:pol method computational aspects}. Techniques for determining the density function of the limiting eigenvalue distribution function and the associated moments are discussed in Sections \ref{sec:level density} and \ref{sec:moments}, respectively. We discuss the relevance of the polynomial method to computational free probability in Section \ref{sec:computational free probability}, provide some applications in Section \ref{sec:pol method examples} and conclude with some open problems in Section \ref{sec:open problems}.

\section{Transform representations}\label{sec:transforms}
We now describe the various ways in which transforms of the empirical distribution function can be encoded and manipulated.
\subsection{The Stieltjes transform and some minor variations}
The Stieltjes transform of the distribution function $F^{A}(x)$ is given by
\begin{equation}\label{eq:mz}
m_{A}(z)=\int \frac{1}{x-z}dF^{A}(x) \qquad \textrm{ for } z\in \mathbb{C}^{+}\setminus \mathbb{R}.
\end{equation}
The Stieltjes transform may be interpreted as the expectation 
\begin{equation*}\label{eq:E mz}
m_{A}(z)=E_{x} \bigg[ \frac{1}{x-z} \bigg],
\end{equation*}
with respect to the random variable $x$ with distribution function $F^{A}(x)$. Consequently, for any invertible function $h(x)$ continuous over the support of $dF^{A}(x)$,
the Stieltjes transform $m_{A}(z)$ can also be written in terms of the distribution of the random variable $y=h(x)$ as
\begin{equation}\label{eq:mz x y}
m_{A}(z)=E_{x} \bigg[ \frac{1}{x-z} \bigg] = E_{y} \bigg[ \frac{1}{h^{\langle -1 \rangle}(y)-z} \bigg],
\end{equation}
where $h^{\langle -1 \rangle}(\cdot)$ is the inverse of $h(\cdot)$ with respect to composition i.e. $h(h^{\langle
-1\rangle}(x))=x$. Equivalently, for $y=h(x)$, we obtain the relationship
\begin{equation}\label{eq:mz y x}
E_{y} \bigg[ \frac{1}{y-z} \bigg]=E_{x} \bigg[ \frac{1}{h(x)-z} \bigg].
\end{equation}
The well-known Stieltjes-Perron inversion formula \cite{akhiezer65a} 
\begin{equation}\label{eq:inversion formula}
f_{A}(x) \equiv dF^{A}(x) =\frac{1}{\pi} \lim _{\xi \rightarrow 0^{+}}{\rm Im} \:m_{A}(x+i\xi).
\end{equation}
can be used to recover the probability density function  $f_{A}(x)$ from the Stieltjes transform. Here and for the remainder of this thesis, the density function is assumed to be distributional derivative of the distribution function. In a portion of the literature on random matrices, the Cauchy transform is defined as
\begin{equation*}
g_{A}(z)=\int \frac{1}{z-x} dF^{A}(x) \qquad for z \in \mathbb{C}^{-1} \setminus \mathbb{R}.
\end{equation*}
The Cauchy transform is related to the Stieltjes transform, as defined in (\ref{eq:mz}), by
\begin{equation}\label{eq:mz gz} g_{A}(z)=-m_{A}(z).
\end{equation}
\subsection{The moment transform} When the probability distribution is compactly supported, the Stieltjes transform can also be expressed as the series expansion
\begin{equation}\label{eq:m moment}
m_{A}(z)=-\dfrac{1}{z}-\sum_{j=1}^{\infty}\dfrac{M_{j}^{A}}{z^{j+1}},
\end{equation}
about $z=\infty$, where $M_{j}^{A} := \int x^{j} dF^{A}(x)$ is the $j$-th moment. The ordinary moment generating function, $\mu_{A}(z)$, is the power series 
\begin{equation}\label{eq:moment transform}
\mu_{A}(z)=\sum_{j=0}^{\infty} M_{j}^{A}\,z^{j},
\end{equation}
with $M_{0}^{A} = 1$. The moment generating function, referred to as the moment transform, is related to the Stieltjes transform by
\begin{equation}\label{eq:mu z}
\mu_{A}(z)=-\frac{1}{z}m_{A}\left(\frac{1}{z}\right).
\end{equation}
The Stieltjes transform can be expressed in terms of the moment transform as
 \begin{equation}\label{eq:m in terms of mu} 
m_{A}(z)=-\frac{1}{z}\mu_{A}\left(\frac{1}{z}\right).
\end{equation}
The eta transform, introduced by Tulino and Verd\`{u} in \cite{tulino04a}, is a minor variation of the moment transform. It can be expressed in terms of the Stieltjes transform as
\begin{equation}
\eta_{A}(z)=\frac{1}{z}m_{A}\left(-\frac{1}{z}\right),
\end{equation}
while the Stieltjes transform can be expressed in terms of the eta transform as
\begin{equation}\label{eq:eta z}
m_{A}(z)=-\frac{1}{z}\eta_{A}\left(-\frac{1}{z}\right).
\end{equation}
\subsection{The R transform}
The R transform is defined in terms of the Cauchy transform as
\begin{equation}
r_{A}(z)=g_{A}^{\langle -1 \rangle}(z)-\frac{1}{z},
\end{equation}
where $g_{A}^{\langle -1 \rangle}(z)$ is the functional inverse of $g_{A}(z)$ with respect to composition. It
will often be more convenient to use the expression for the R transform in terms of the Cauchy transform given
by
\begin{equation}\label{eq:r transform}
r_{A}(g)=z(g)-\frac{1}{g}.
\end{equation}
The R transform can be written as a power series whose coefficients $K_{j}^{A}$ are known as the ``free
cumulants.'' For a combinatorial interpretation of free cumulants, see \cite{speicher97a}. Thus the R transform is the (ordinary) free cumulant generating function
\begin{equation}
r_{A}(g)=\sum_{j=0}^{\infty}\: K_{j+1}^{A}\,g^{j}.
\end{equation}
\subsection{The S transform}
 The S transform is relatively more complicated. It is defined as
 \begin{equation}\label{eq:s1}
s_{A}(z) = \frac{1+z}{z}\, \Upsilon_{A}^{\langle -1 \rangle}(z)
 \end{equation}
 where $\Upsilon_{A}(z)$ can be written in terms of the Stieltjes transform $m_{A}(z)$ as
 \begin{equation}\label{eq:s2}
\Upsilon_{A}(z)=-\frac{1}{z}m_{A}(1/z)-1.
 \end{equation}
This definition is quite cumbersome to work with because of the functional inverse in (\ref{eq:s1}). It also
places a technical restriction (to enable series inversion) that $M_{1}^{A} \neq 0$. We can, however, avoid this
by expressing the S transform algebraically in terms of the Stieltjes transform as shown next.
We first plug in $\Upsilon_{A}(z)$ into the left-hand side of (\ref{eq:s1}) to obtain
\begin{equation*}
s_{A}(\Upsilon_{A}(z)) = \frac{1+\Upsilon_{A}(z)}{\Upsilon_{A}(z)}\, z.
\end{equation*}
This can be rewritten in terms of $m_{A}(z)$ using the relationship in (\ref{eq:s2}) to obtain
\begin{equation*}
s_{A}(-\frac{1}{z}m(1/z)-1)=\frac{z\,m(1/z)}{m(1/z)+z}
\end{equation*}
or, equivalently:
\begin{equation}\label{eq:s3}
s_{A}(-z\,m(z)-1)=\frac{m(z)}{z\,m(z)+1}.
\end{equation}
We now define $y(z)$ in terms of the Stieltjes transform as $y(z)=-z\,m(z)-1$. It is clear that $y(z)$ is an
invertible function of $m(z)$. The right hand side of (\ref{eq:s3}) can be rewritten in terms of $y(z)$ as
\begin{equation}\label{eq:s4}
s_{A}(y(z)) = -\frac{m(z)}{y(z)}=\frac{m(z)}{z\,m(z)+1}.
\end{equation}
Equation (\ref{eq:s4}) can be rewritten to obtain a simple relationship between the Stieltjes transform and the
S transform
\begin{equation}\label{eq:mz and sy}
m_{A}(z) = -y\,s_{A}(y).
\end{equation}
Noting that $y=-z\,m(z)-1$ and $m(z) = -y\,s_{A}(y)$ we obtain the relationship
\begin{equation*}
y = z\,y\,s_{A}(y)-1
\end{equation*}
or, equivalently
\begin{equation}\label{eq:s5}
z=\frac{y+1}{y\,s_{A}(y)}.
\end{equation}
\section{Algebraic distributions}\label{sec:alg distributions}
\begin{notation}[Bivariate polynomial]
Let $\luvs{}$ denote a bivariate polynomial of degree $\Du$ in $u$ and $\Dv$ in $v$ defined as 
\begin{equation}\label{eq:bivariate polynomial}
\luvs{} \equiv \luvs{}(\cdot,\cdot) = \sum_{j=0}^{\Du}\sum_{k=0}^{\Dv} c_{jk}\, u^{j}\,v^k = \sum_{j=0}^{\Du}\, l_{j}(v)\, u^j.
\end{equation}
The scalar coefficients $c_{jk}$ are real valued. \\
\end{notation}
The two letter subscripts for the bivariate polynomial $\luvs{}$ provide us with  a convention of which dummy variables we will use. We will generically use the first letter in the subscript to represent a transform of the density with the second letter acting as a mnemonic for the dummy variable associated with the transform. By consistently using the same pair of letters to denote the bivariate polynomial that encodes the transform and the associated dummy variable, this abuse of notation allows us to readily identify the encoding of the distribution that is being manipulated.
\begin{remark}[Irreducibility]
Unless otherwise stated it will be understood that $\luvs{}(u,v)$ is ``irreducible'' in the sense that the conditions:\\
\begin{itemize}
\item $l_{0}(v),\ldots, l_{\Du}(v)$ have no common factor involving $v$,
\item $l_{\Du}(v) \neq 0$,
\item ${\rm disc}_{{\rm L}}(v) \neq 0$,
\end{itemize}
are satisfied, where ${\rm disc}_{L}(v)$ is the discriminant of $\luv{}$ thought of as a polynomial in $v$.\\
\end{remark}
\noindent We are particularly focused on the solution ``curves,'' $u_{1}(v),\ldots,u_{\Du}(v)$, \ie,
\[
L_{\textrm{uv}}(u,v)=l_{\Du}(v)\, \prod_{i=1}^{D_{{\rm u}}}\left(u-u_{i}(v)\right).
\]
Informally speaking, when we refer to the bivariate polynomial equation $\luv{}=0$ with solutions $u_{i}(v)$ we are actually considering the equivalence class of rational functions with this set of solution curves. \\
\begin{remark}[Equivalence class]\label{remark:equivalence class}
The equivalence class of $\luv{}$ may be characterized as functions of the form $\luv{}g(v)/h(u,v)$ where $h$ is relatively prime to $\luv{}$ and $g(v)$ is not identically $0$.\\ 
\end{remark}
\noindent A few technicalities (such as poles and singular points) that will be catalogued later in Section \ref{sec:pol method computational aspects} remain, but this is sufficient for allowing us to introduce rational transformations of the arguments and continue to use the language of polynomials.
\begin{definition}[Algebraic distributions]
Let $F(x)$ be a probability distribution function and $f(x)$ be its distributional derivative (here and henceforth). Consider the Stieltjes transform $m(z)$ of the distribution function, defined as
\begin{equation}
m(z) = \int \dfrac{1}{x-z}dF(x) \qquad \textrm{for } z \in \mathbb{C}^{+} \setminus \mathbb{R}.
\end{equation}
If there exists a bivariate polynomial $\lmzs{}$ such that $\lmzs{}(m(z),z)=0$ then we refer to $F(x)$ as algebraic (probability) distribution function, $f(x)$ as an algebraic (probability) density function and say the $f \in \mathcal{P}_{\textrm{alg}}$. Here $\Palg$ denotes the class of algebraic (probability) distributions. 
\end{definition}
\begin{definition}[Atomic distribution]\label{definition:atomic density}
Let $F(x)$ be a probability distribution function of the form
\[
F(x) = \sum_{i=1}^{K} p_{i}\,\mathbb{I}_{[\lambda_{i},\infty)},
\]
where the $K$ atoms at $\lambda_{i} \in \mathbb{R}$ have (non-negative) weights $p_{i}$  subject to $\sum_{i} p_{i}=1$ and $\mathbb{I}_{[x,\infty)}$ is the indicator (or characteristic) function of the set $[x,\infty)$. We refer to $F(x)$ as an atomic (probability) distribution function. Denoting its distributional derivative by $f(x)$, we say that $f(x) \in \Patom$. Here $\Patom$ denotes the class of atomic distributions.
\end{definition}
%\begin{definition}[Positive semi-definite algebraic density]
%Algebraic densities supported over the positive real axis are referred to as positive semi-definite. We denote the class of positive %semi-definite densities by $\PalgP$. Thus, if  $\roc{A} \subset [0,\infty)$ then we say that $f_{A}(x) \in \PalgP$. \\
%\end{definition}
\begin{example}
An atomic probability distribution, as in Definition \ref{definition:atomic density}, has a Stieltjes transform 
\[
m(z) = \sum_{i=1}^{K}\dfrac{p_{i}}{\lambda_{i}-z}
\]
which is the solution of the equation $\lmz{}=0$ where $$\lmz{} \equiv \prod_{i=1}^{K} (\lambda_{i}-z)\, m -  \sum_{i=1}^{K} \prod_{\begin{subarray}{l}j \neq i\\ j =1\end{subarray}}^{K} p_{i}(\lambda_{j}-z).$$
Hence it is an algebraic distribution; consequently $\Patom \subset \Palg$.
\end{example}
\begin{example}
The Cauchy distribution whose density
\begin{equation*}
f(x) = \dfrac{1}{\pi(x^2+1)},
\end{equation*}
has a Stieltjes transform $m(z)$ which is the solution of the equation $\lmz{}=0$ where  $$\lmz{} \equiv \left( {z}^{2}+1 \right) {m}^{2}+2\,z\,m+1.$$ Hence it is an algebraic distribution.
\end{example}
It is often the case that the probability density functions of algebraic distributions, according to our definition, will also be algebraic functions themselves. We conjecture that this is a necessary but not sufficient condition. We show that it is not sufficient by providing the counter-example below.\\
\begin{counterexample}
Consider the quarter-circle distribution with density function 
\begin{equation*}
f(x) = \dfrac{\sqrt{4-x^2}}{\pi} \qquad \textrm{ for } x \in [0,2].
\end{equation*}
Its Stieltjes transform :
\begin{equation*} 
m(z) =  -\dfrac{4-2\,\sqrt {-{z}^{2}+4}\ln  \left( -{\frac {2+\sqrt {-{z}^{2}+4}}{z}} \right) +z\pi }{2\pi},
\end{equation*}
is clearly not an algebraic function. Thus $f(x) \notin \Palg$. \\
\end{counterexample}
\subsection{Implicit representations of algebraic distributions}
\begin{figure}[t]
\centering
\includegraphics[width=5.7in]{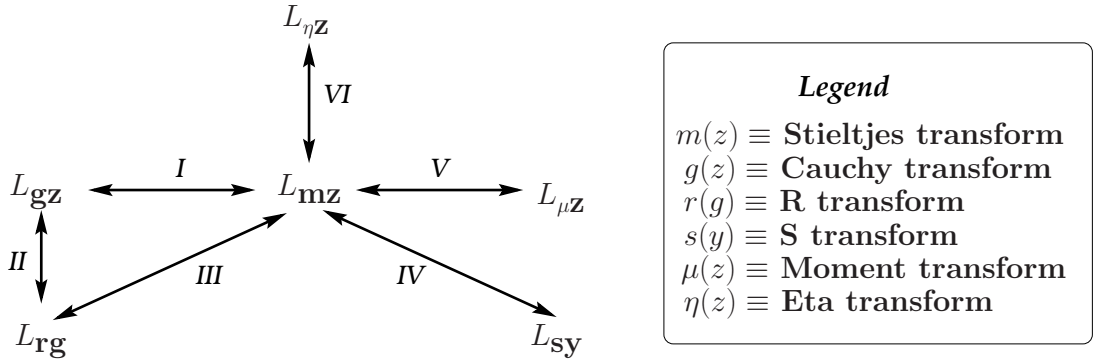}
\caption{The six interconnected bivariate polynomials; transformations between the polynomials, indicated by
the labelled arrows, are given in Table \ref{tab:transforms}.} 
\label{fig:all transforms}
\end{figure}

We now define six interconnected bivariate polynomials denoted by $\lmzs{}$, $\lgzs{}$, $\lrgs{}$, $\lsys{}$,
$\lmuzs{}$, and $\letazs{}$.  We assume that $\luv{}$ is an irreducible bivariate polynomial of the form in (\ref{eq:bivariate polynomial}).
The main protagonist of the transformations we consider is the bivariate polynomial $\lmzs{}$ which implicitly defines the Stieltjes transform $m(z)$ via the equation $\lmzs{}(m,z)=0$. Starting off with this polynomial we can obtain the polynomial $\lgzs{}$ using the relationship in (\ref{eq:mz gz}) as
\begin{equation}\label{eq:Lgz from Lmz}
\lgz{} =\lmzs{}(-g,z).
\end{equation}
Perhaps we should explain our abuse of notation once again, for the sake of clarity. Given any one polynomial, all the other polynomials can be obtained. The
two letter subscripts not only tell us which of the six polynomials we are focusing on, it provides a convention
of which dummy variables we will use. The first letter in the subscript represents the transform; the second letter is a mnemonic for the variable associated with the transform that we use consistently in the software based on this framework.
With this notation in mind, we can obtain the polynomial $\lrgs{}$ from $\lgzs{}$ using (\ref{eq:r transform}) as
\begin{equation}\label{eq:Lrg from Lgz}
\lrg{} =\lgzs{}\left(g,r+\dfrac{1}{g}\right).
\end{equation}
Similarly, we can obtain the bivariate polynomial $\lsys{}$ from $\lmzs{}$ using the expressions in (\ref{eq:mz and sy}) and (\ref{eq:s5}) to obtain the relationship
\begin{equation}\label{eq:Lsy from Lmz}
\lsys{}=\lmzs{}\left(-y\,s,\frac{y+1}{sy}\right).
\end{equation}
Based on the transforms discussed in Sectin \ref{sec:transforms}, we can derive transformations between
additional pairs of bivariate polynomials represented by the bidirectional arrows in Figure \ref{fig:all
transforms} and listed in the third column of Table \ref{tab:transforms}. Specifically, the expressions in
(\ref{eq:mu z}) and (\ref{eq:eta z}) can be used to derive the transformations between $\lmzs{}$ and $\lmuzs{}$
and $\lmzs{}$ and $\letazs{}$ respectively. 
The fourth column of Table \ref{tab:transforms} lists the \matlab function, implemented using its \maple based Symbolic Toolbox, corresponding to the bivariate polynomial transformations represented in Figure \ref{fig:all transforms}. In the \matlab functions, the function \texttt{irreducLuv(u,v)} listed in Table \ref{tab:irreducible code} ensures that the resulting bivariate polynomial is irreducible by clearing the denominator and making the resulting polynomial square free. \\

\vspace{0.15in}
\noindent\textbf{Example}: Consider an atomic probability distribution with
\begin{equation}\label{eq:edf example}
F(x)=0.5 \, \mathbb{I}_{[0,\infty)} + 0.5 \,\mathbb{I}_{[1,\infty)},
\end{equation}
whose Stieltjes transform
\begin{equation*}
m(z) =\dfrac{0.5}{0-z} + \dfrac{0.5}{1-z},
\end{equation*}
is the solution of the equation
\begin{equation*}
m(0-z)(1-z)-0.5(1-2z)=0,
\end{equation*}
or equivalently, the solution of the equation $\lmz{} = 0$ where 
\begin{equation}\label{eq:LmzA eg}
\lmz{}\equiv m(2\,z^2-2\,z)-(1-2z).
\end{equation}
We can obtain the bivariate polynomial $\lgz{}$ by applying the transformation in (\ref{eq:Lgz from Lmz}) to the bivariate polynomial $\lmzs{}$ given by (\ref{eq:LmzA eg}) so that
\begin{equation}\label{eq:LgzA eg}
\lgz{}= -g(2\,z^2-2\,z)-(1-2z).
\end{equation}
Similarly, by applying the transformation in (\ref{eq:Lrg from Lgz})  we obtain
\begin{equation}\label{eq:example Lrg rational}
\lrg{} = -g\left(2\left(r+\frac{1}{g}\right)-2\left(r+\frac{1}{g}\right)^2\right)-\left(1-2\left(r+\frac{1}{g}\right)\right).
\end{equation}
which, on clearing the denominator and invoking the equivalence class representation of our polynomials (see Remark \ref{remark:equivalence class}), gives us the irreducible bivariate polynomial
\begin{equation}
\lrg{} = -1+2\,g{r}^{2}+ \left( 2-2\,g \right) r.
\end{equation}
\vspace{0.25cm}
By applying the transformation in (\ref{eq:Lsy from Lmz}) to the bivariate polynomial $\lmzs{}$, we obtain 
\begin{equation*}
\lsys{} \equiv (-s\,y)\left(2\,\frac{y+1}{sy}-2\,\left(\frac{y+1}{sy}\right)^2\right)-\left(1-2\frac{y+1}{sy}\right)
\end{equation*}
which on clearing the denominator gives us the irreducible bivariate polynomial
\begin{equation}
\lsy{A} = \left( 1+2\,y \right) s-2-2\,y.
\end{equation}
Table \ref{tab:couple of examples} tabulates the six bivariate polynomial encodings in Figure \ref{fig:all transforms} for the distribution in (\ref{eq:edf example}), the semi-circle distribution for Wigner matrices and the Mar\v{c}enko-Pastur distribution for Wishart matrices.
\begin{table}[t]
\centering
\begin{tabular}{|c|l|}
\hline
Procedure & \matlab Code \\
\hline
                         & {\small \texttt{function Luv = irreducLuv(Luv,u,v)}} \\
                         & \\ 
Simplify and clear the denominator  & {\small \texttt{L = numden(simplify(expand(Luv)));}}\\
                                      & {\small \texttt{L = Luv / maple('gcd',L,diff(L,u));}} \\
     Make square free                  & {\small \texttt{L = simplify(expand(L));}}\\
                                       & {\small \texttt{L = Luv / maple('gcd',L,diff(L,v));}} \\
Simplify                                    & {\small \texttt{Luv = simplify(expand(L));}} \\  
\hline
\end{tabular}
\caption{Making $\luvs{}$ irreducible.}
\label{tab:irreducible code}
\end{table}
\begin{table}[b]
\centering
\subtable[The atomic distribution in (\ref{eq:edf example}).]{\label{tab:example 1 summary}
\begin{tabular}{|l||l|}
\hline
$L$ & Bivariate Polynomials \\
\hline
$\lmzs{}$ & $m(2\,z^2-2\,z)-(1-2z)$ \\
$\lgzs{}$ & $-g(2\,z^2-2\,z)-(1-2z)$\\
$\lrgs{}$ & $-1+2\,g{r}^{2}+ \left( 2-2\,g \right) r$ \\
$\lsys{}$ & $ \left( 1+2\,y \right) s-2-2\,y$ \\
$\lmuzs{}$ &$ \left( -2+2\,z \right)\mu+2-z$ \\
$\letazs{}$ & $ \left( 2\,z+2 \right) \eta-2-z$ \\
\hline
\end{tabular}
} 
\hspace{0.40in}
\subtable[The Mar\v{c}enko-Pastur distribution.]{\label{tab:Wishart bivariate}
\begin{tabular}{|l||l|}
\hline
$L$ & Bivariate Polynomials \\
\hline
$\lmzs{}$ & $cz{m}^{2}- \left( 1-c-z \right) m+1$ \\
$\lgzs{}$ & $cz{g}^{2}+\left( 1-c-z \right) g+1$ \\
$\lrgs{}$ & $\left( cg-1 \right) r+1$\\
$\lsys{}$ & $\left( cy+1 \right) s-1$\\
$\lmuzs{}$ &$\mu^{2}zc- \left( zc+1-z \right)\mu+1$\\
$\letazs{}$ & ${\eta}^{2}zc+ \left(-zc+1-z \right) \eta-1$ \\
\hline
\end{tabular}
} \\[0.3in]
\hspace{0.5in}
\subtable[The semi-circle distribution.]{\label{tab:Hermite bivariate}
\begin{tabular}{|l||l|}
\hline
$L$ & Bivariate polynomials  \\[0.1cm]
\hline \hline
$\lmzs{}$ & $m^{2}+m\,z+1$\\
$\lgzs{}$ & $g^{2}-g\,z+1$\\
$\lrgs{}$ & $r-g$ \\
$\lsys{}$ & $s^{2}\,y-1$\\
$\lmuzs{}$ & $\mu^{2}{z}^{2}-\mu+1$\\
$\letazs{}$ & ${z}^{2}{\eta}^{2}-\eta+1$ \\
\hline
\end{tabular}
}
\caption{Bivariate polynomial representations of some algebraic distributions.}
\label{tab:couple of examples}
\end{table}
\begin{table}[p]
\centering
\captionstyle{flushleft}
\includegraphics[width=5.3in]{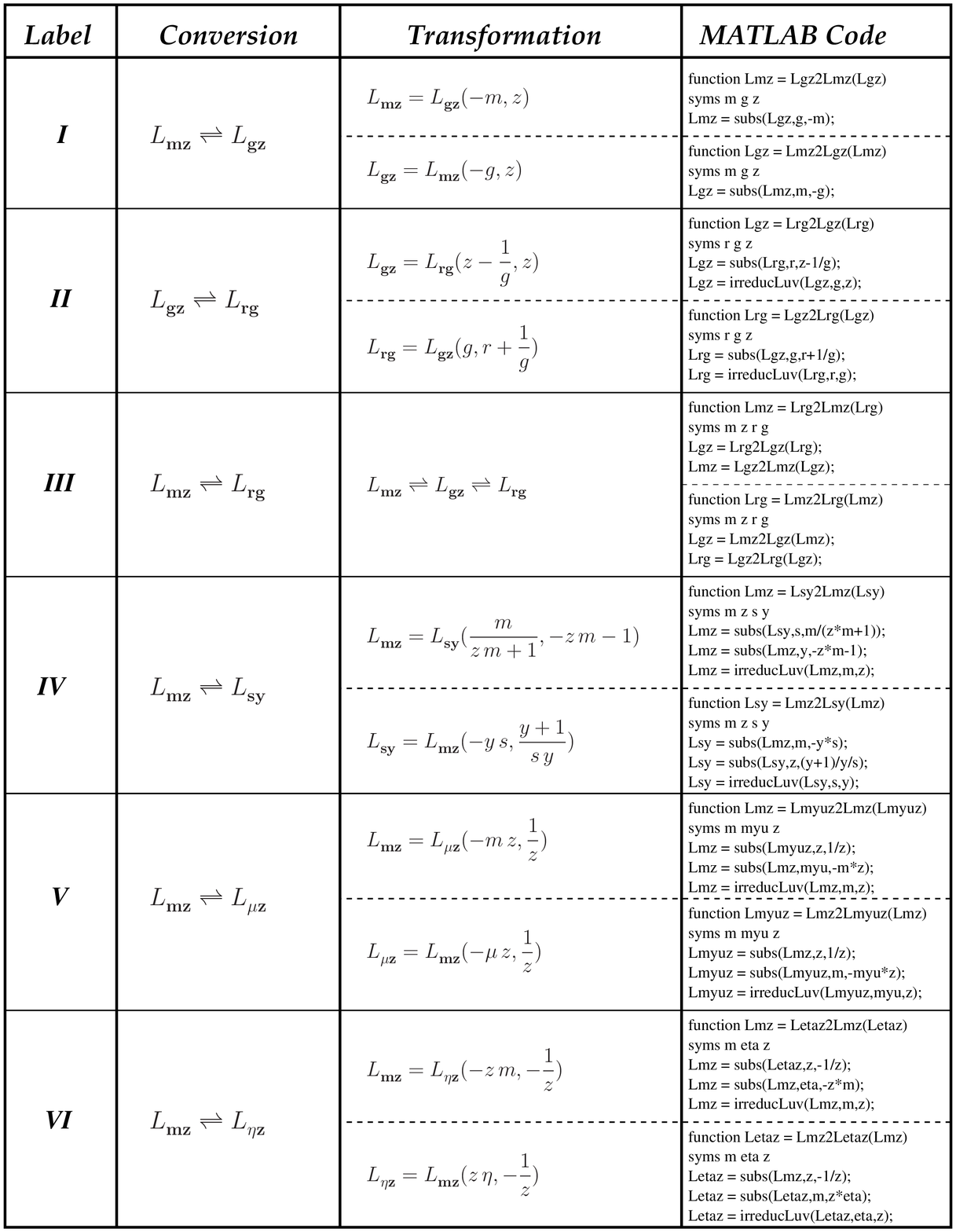}
\caption{Transformations between the different bivariate polynomials. As a guide to \MATLAB notation, the command
\texttt{syms} declares a variable to be symbolic while the command \texttt{subs} symbolically substitutes every
occurrence of the second argument in the first argument with the third argument. Thus, for example, the command
\texttt{y=subs(x-a,a,10)} will yield the output \texttt{y=x-10} if we have previously declared $x$ and $a$ to be
symbolic using the command \texttt{syms x a}.} \label{tab:transforms}
\end{table}
\clearpage
\section{Algebraic operations on algebraic functions}\label{sec:kronecker}
Algebraic functions are closed under addition and multiplication. Hence we can add (or multiply) two algebraic functions and obtain another algebraic function. We show, using purely matrix theoretic arguments, how to  obtain the polynomial equation  whose solution is the sum (or product) of two algebraic functions without ever actually computing the individual functions. In  Section \ref{sec:resultant connection}, we interpret this computation using the concept of resultants \cite{sturmfels97a} from elimination theory.  These tools will feature prominently in Section \ref{sec:algebraic random matrices} when we encode the transformations of the random matrices as algebraic operations on the appropriate form of the  bivariate polynomial that encodes their limiting eigenvalue distributions. \\
\subsection{Companion matrix based computation}
\begin{definition}[Companion Matrix]
The companion matrix ${\bf C}_{a(x)}$ to a monic polynomial
\begin{equation*}
a(x) \equiv a_0 + a_1\,x+ \ldots + a_{n-1}\,x^{n-1} + x^{n}
\end{equation*}
is the $n \times n$ square matrix
\begin{small}
\begin{equation*}
{\bf C}_{a(x)} =
\begin{bmatrix}
0 &\ldots & \ldots & \ldots & -a_{0} \\
1 & \cdots & \cdots & \cdots & -a_{1} \\
0 & \ddots &        &        & -a_{2} \\
\vdots &         &\ddots  &        & \vdots \\
0       & \ldots & \ldots  & 1     & -a_{n-1}
\end{bmatrix}
\end{equation*}
\end{small}
with ones on the sub-diagonal and the last column given by the negative coefficients of $a(x)$.\\
\end{definition}
\begin{remark}
The eigenvalues of the companion matrix are the solutions of the equation $a(x)=0$. This is intimately related to the observation that the characteristic polynomial of the companion matrix equals $a(x)$, i.e.,
\begin{equation*}
a(x) = \det(x\,{\bf I}_{n}-{\bf C}_{a(x)}). \\
\end{equation*}
\end{remark}
\noindent Consider the bivariate polynomial $\luvs{}$ as in (\ref{eq:bivariate polynomial}). By treating it as a polynomial in $u$ whose coefficients are polynomials in $v$, \ie, by rewriting it as
\begin{equation}\label{eq:Cuv Luv}
\luv{}\equiv \sum_{j=0}^{\Du} l_{j}(v)\,u^{j},
\end{equation}
we can create a companion matrix $\Cuv{u}$  whose characteristic polynomial as a function of $u$ is the
bivariate polynomial $\luvs{}$. The companion matrix $\Cuv{u}$ is the $\Du \times \Du$ matrix in Table
\ref{tab:companion matrix}.
\begin{table}[h]
\centering
\begin{tabular}{|c||c|}
\hline $\Cuv{u}$ & \MATLAB code \\
\hline
\begin{small}
\begin{tabular}{c}
$\left[
\begin{array}{ccccc}
0 &\ldots & \ldots & \ldots & -l_{0}(v)/l_{\Du}(v) \\
1 & \cdots & \cdots & \cdots & - l_{1}(v)/l_{\Du}(v)\\
0 & \ddots &        &        & -l_{2}(v)/l_{\Du}(v) \\
\vdots &         &\ddots  &        & \vdots \\
0       & \ldots & \ldots  & 1 & -l_{\Du\!-1}(v)/l_{\Du}(v)
\end{array}\right]
$
\end{tabular}
\end{small}
&
\begin{footnotesize}
\begin{tabular}{l}
function Cu = Luv2Cu(Luv,u)\\
Du = double(maple('degree',Luv,u));\\
LDu = maple('coeff',Luv,u,Du);\\
Cu = sym(zeros(Du))+ ..\\
\phantom{iiiiiii}+diag(ones(Du-1,1),-1));\\
for Di = 0:Du-1\\
\phantom{iiiiii}LtuDi = maple('coeff',Lt,u,Di);\\
\phantom{iiiiii}Cu(Di+1,Du) = -LtuDi/LDu;\\
end\\
\end{tabular}
\end{footnotesize}\\
\hline
\end{tabular}
\caption{The companion matrix $\Cuv{u}$, with respect to $u$, of the bivariate polynomial $\luvs{}$ given by (\ref{eq:Cuv
Luv}).} \label{tab:companion matrix}
\end{table}
\begin{remark}
Analogous to the univariate case, the characteristic polynomial  of $\Cuv{u}$ is $\det(u\,{\bf I}-\Cuv{u}) = \luv{}/l_{\Du}(v)^{\Du}$. Since $l_{\Du}(v)$ is not identically zero, we say that $\det(u\,{\bf I}-\Cuv{u}) = \luv{}$ where the equality is understood to be with respect to the equivalence class of $\luvs{}$ as in Remark \ref{remark:equivalence class}.  The eigenvalues of $\Cuv{u}$ are the solutions of the algebraic equation $\luv{}=0$; specifically, we obtain the algebraic function $u(v)$. \\
\end{remark}
\begin{definition}[Kronecker product]
If ${\bf A}_{m}$ (with entries $a_{ij}$) is an $m \times m$ matrix and ${\bf B}_{n}$ is an $n \times n$ matrix then the Kronecker (or tensor) product of ${\bf A}_{m}$ and ${\bf B}_{n}$, denoted by ${\bf A}_{m} \otimes {\bf B}_{n}$, is the $mn \times mn$ matrix defined as:
\begin{small}
\begin{equation*}
{\bf A}_{m} \otimes {\bf B}_{n} =
\begin{bmatrix}
a_{11}{\bf B}_{n} & \ldots &a_{1n}{\bf B}_{n} \\
\vdots & \ddots & \vdots \\
a_{m1}{\bf B}_{n} & \ldots & a_{mn}{\bf B}_{n}
\end{bmatrix}\\
\end{equation*}
\end{small}
\end{definition}
\begin{lemma}\label{lemma: tensor product eigenvalues}
 If $\alpha_{i}$ and $\beta_{j}$ are the  eigenvalues of ${\bf A}_{m}$ and ${\bf B}_{n}$ respectively, then
\begin{enumerate}
\item $\alpha_{i}+\beta_{j}$ is an eigenvalue of $({\bf A}_{m}\otimes {\bf I}_{n})+({\bf I}_{m}\otimes {\bf B}_{n})$,
\item $\alpha_{i}\,\beta_{j}$ is an eigenvalue of ${\bf A}_{m}\otimes {\bf B}_{n}$,
\end{enumerate}
for $i=1,\ldots,m$, $j=1,\ldots,n$. 
\end{lemma}
\begin{proof}
The statements are proved in \cite[Theorem 4.4.5]{horn91a} and \cite[Theorem 4.2.12]{horn91a}.
\end{proof}
\begin{proposition}\label{prop:addition of algebraic}
Let $u_{1}(v)$ be a solution of the algebraic equation $\luv{1}=0$, or equivalently an eigenvalue of the $\Dus{1}
\times \Dus{1}$ companion matrix $\Cuv{u_{1}}$. Let $u_{2}(v)$ be a solution of the algebraic equation $\luv{2} =0 $,
or equivalently an eigenvalue of the $\Dus{2} \times \Dus{2}$ companion matrix $\Cuv{u_{2}}$. Then
\begin{enumerate}
\item $u_{3}(v)=u_{1}(v)+u_{2}(v)$ is an eigenvalue of the matrix $\Cuv{u_3}=\left(\Cuv{u_{1}} \otimes
{\bf I}_{\Dus{2}}\right)+\left({\bf I}_{\Dus{1}} \otimes \Cuv{u_{2}}\right)$,
 \item $u_{3}(v)=u_{1}(v)u_{2}(v)$ is an eigenvalue of the matrix $\Cuv{u_3}=\Cuv{u_{1}}\otimes
 \Cuv{u_{2}}$.
\end{enumerate}
Equivalently $u_3(v)$ is a solution of the algebraic equation $\luvs{3} =0$ where $ \luvs{3}=\det(u\,{\bf I}-\Cuv{u_3})$.
\end{proposition}
\begin{proof}
This follows directly from Lemma \ref{lemma: tensor product eigenvalues}.
\end{proof}
\vspace{0.25cm}
We represent the binary addition and multiplication operators on the space of algebraic functions by the symbols $\bp{u}$ and $\bt{u}$ respectively. We define addition and multiplication as in Table \ref{tab:bp and bm operators} by applying Proposition \ref{prop:addition of algebraic}.  Note that the subscript `u' in $\bp{u}$ and $\bt{u}$ provides us with an indispensable convention of which dummy variable we are using. Table \ref{tab:example of Luv operations} illustrates the $\bp{}$ and $\bt{}$ operations on a pair of bivariate polynomials and underscores the importance of the symbolic software developed. The $(\Du+1) \times (\Dv+1)$ matrix ${\bf T}_{{\rm uv}}$ lists only the coefficients $c_{ij}$ for the term $u^{i}\,v^{j}$ in the polynomial $\luv{}$. Note that the indexing for $i$ and $j$ starts with zero.
\begin{table}[h]
\centering
\includegraphics[width=5.75in]{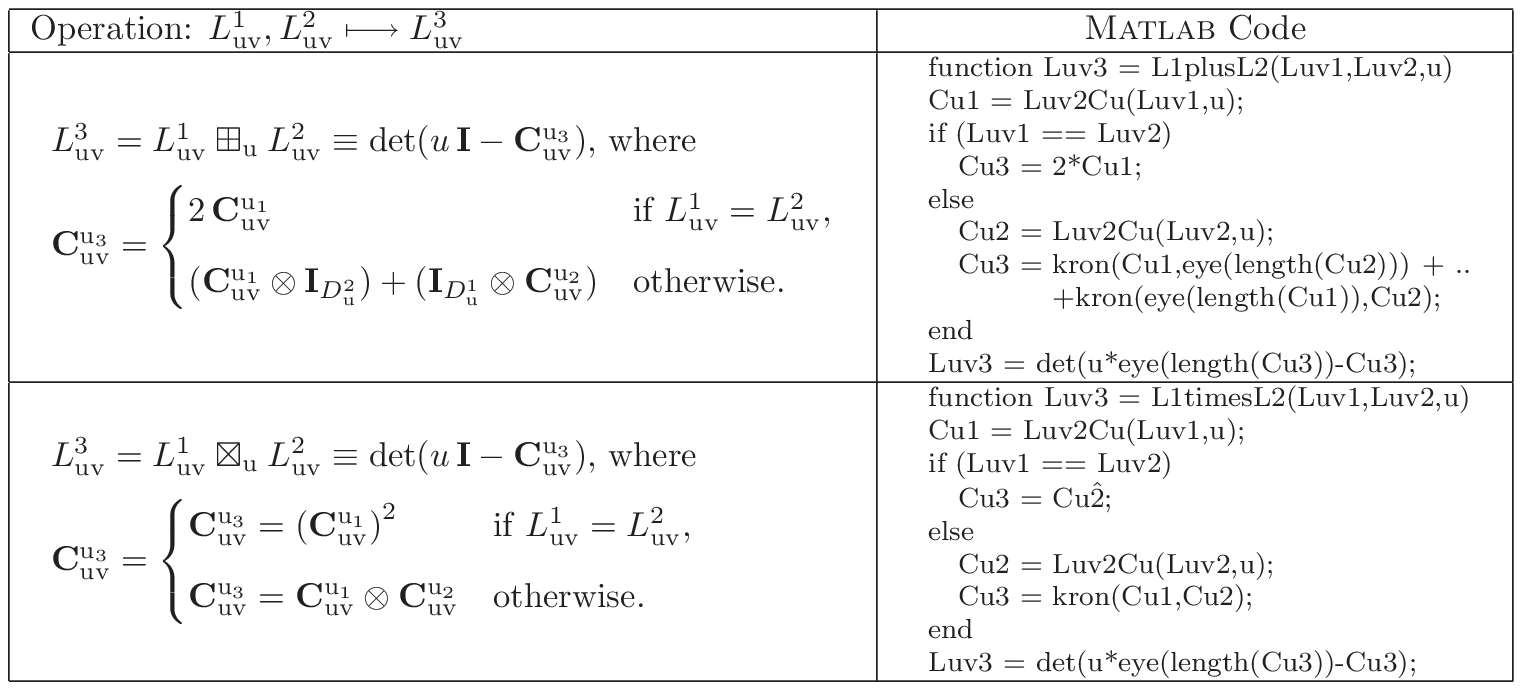}
\caption{Formal and computational description of the $\bp{u}$ and $\bt{u}$ operators acting on the bivariate
polynomials $\luv{1}$ and $\luv{2}$ where $\Cuv{u_1}$ and $\Cuv{u_2}$ are their corresponding companion matrices
constructed as in Table \ref{tab:companion matrix} and $\otimes$ is the matrix Kronecker product.} \label{tab:bp
and bm operators}
\end{table}
\begin{table}
\centering
\includegraphics[width=5.75in]{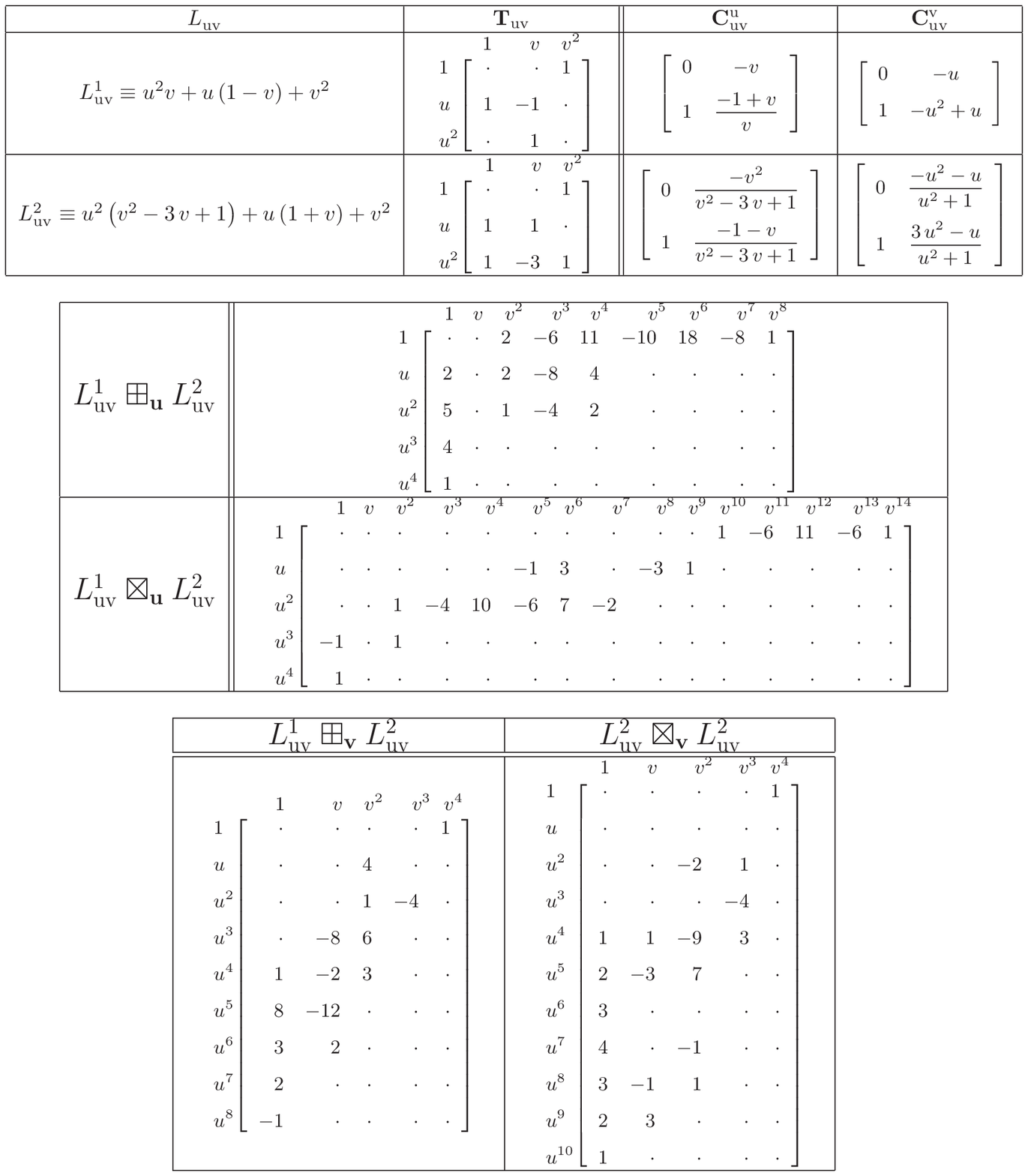}
\caption{Examples of $\bp{}$ and $\bt{}$ operations on a pair of bivariate polynomials, $\luvs{1}$ and $\luvs{2}$.}
\label{tab:example of Luv operations}
\end{table} 
\subsection{Resultants based computation} \label{sec:resultant connection}
 Addition (and multiplication) of algebraic functions produces another algebraic function. We now demonstrate how the concept of resultants from elimination theory can be used to obtain the polynomial whose zero set is the required algebraic function.
\begin{definition}[Resultant]\label{def:resultant}
Given a polynomial
\[
a(x) \equiv a_{0} + a_{1}\,x + \ldots + a_{n-1}\,x^{n-1}+a_{n}x^{n}
\]
of degree $n$ with roots $\alpha_{i}$, for $i=1,\ldots,n$ and a polynomial
\[
b(x) \equiv b_{0} + b_{1}\,x + \ldots + b_{m-1}\,x^{m-1}+b_{m}x^{m}
\]
of degree $m$ with roots $\beta_{j}$, for $j=1,\ldots,m$, the resultant is defined as
 \[
{\rm Res}_{\,{\rm x}}\left(a(x)\,,\,b(x)\right)= a_{n}^{m}\,b_{m}^{n} \prod_{i=1}^{n} \prod_{j=1}^{m}(\beta_{j}-\alpha_{i}).
\]
\end{definition}
From a computational standpoint, the resultant can be directly computed from the coefficients of the polynomials itself. The computation involves the formation of the Sylvester matrix and exploiting an identity that relates the determinant of the Sylvester matrix to the resultant.
\begin{definition}[Sylvester matrix]
Given polynomials $a(x)$ and $b(x)$ with degree $n$ and $m$ respectively and coefficients as in Definition \ref{def:resultant}, the Sylvester matrix is the $(n+m) \times (n+m)$ matrix 
\begin{equation*}
{\bf S}(a,b) = \begin{bmatrix}
a_{n} & 0 & \cdots  &0 & 0 & b_{m} & 0 &\cdots & 0 & 0 \\
a_{n-1} & a_{n} & \cdots  &0 & 0 & b_{m-1} & b_{m} &\cdots & 0 & 0 \\
\hdots & \hdots & \cdots & \hdots & \hdots & \hdots & \hdots & \cdots & \hdots & \hdots \\
0 & 0 & \cdots & a_{0} & a_{1} & 0 & 0 & \cdots & b_{0} & b_{1} \\
0 & 0 & \cdots & 0 & a_{0} & 0 & 0 & \cdots & 0 & b_{0} \\ 
\end{bmatrix}
\end{equation*}\\
\end{definition}
\begin{proposition}\label{prop:resultant identity}
The resultant of two polynomials $a(x)$ and $b(x)$ is related to the determinant of the Sylvester matrix by
\begin{equation*}
\det({\bf S}(a,b)) = {\rm Res}_{\,{\rm x}}\left(a(x)\,,\,b(x)\right)
\end{equation*}
\end{proposition}
\begin{proof}
This identity can be proved using standard linear algebra arguments. A proof may be found in \cite{akritas93a}.
\end{proof}
For our purpose, the utility of this definition is that the $\bp{u}$ and $\bt{u}$ operations can be expressed in terms of resultants. Suppose we are given two bivariate polynomials $\luvs{1}$ and $\luvs{2}$. By using the definition of the resultant and treating the bivariate polynomials as polynomials in $u$ whose coefficients are polynomials in $v$, we obtain the identities
\begin{empheq}[
box=\setlength{\fboxrule}{1pt}\fbox]{equation}\label{eq:resultant add}
\luvs{3}(t,v)=\luvs{1} \bp{u} \luvs{2} \equiv {\rm Res}_{\,{\rm u}} \left(\luvs{1}(t-u,v)\,{\bf ,}\,\luvs{2}(u,v)\right),
\end{empheq}
\noindent and
\begin{empheq}[
box=\setlength{\fboxrule}{1pt}\fbox]{equation}\label{eq:resultant product}
\luvs{3}(t,v)=\luvs{1} \bt{u} \luvs{2} \equiv {\rm Res}_{\,{\rm u}} \left(u^{\Dus{1}}\luvs{1}(t/u,v)\,{\bf,}\,\luvs{2}(u,v)\right),
\end{empheq}
where $\Dus{1}$ is the degree of $\luvs{1}$ with respect to $u$. By Proposition \ref{prop:resultant identity}, evaluating the $\bp{u}$ and $\bt{u}$ operations via the resultant formulation involves computing the determinant of the $(\Dus{1} + \Dus{2}) \times (\Dus{1} + \Dus{2})$ Sylvester matrix. When $\luvs{1} \neq \luvs{2}$, this results in a steep computational saving relative to the companion matrix based formulation in Table \ref{tab:bp and bm operators} which involves computing the determinant of a $(\Dus{1}\Dus{2}) \times (\Dus{1}\Dus{2})$ matrix. Fast algorithms for computing the resultant exploit this and other properties of the Sylvester matrix formulation. In \maple, the computation  $\luvs{3} = \luvs{1} \bp{u} \luvs{2}$ may be performed using the command:
\begin{small}
\begin{verbatim}
   Luv3 = subs(t=u,resultant(subs(u=t-u,Luv1),Luv2,u));
\end{verbatim}
\end{small}
\noindent The computation $\luvs{3} = \luvs{1} \bt{u} \luvs{2}$ can be performed via the sequence of commands:
\begin{small}
\begin{verbatim}
   Du1 = degree(Luv1,u); 
   Luv3 =  subs(t=u,resultant(simplify(u^Du1*subs(u=t/u,Luv1)),Luv2,u));
\end{verbatim}
\end{small}
When $\luvs{1} = \luvs{2}$, however, the $\bp{u}$ and $\bt{u}$ operations are best performed using the companion matrix formulation in Table \ref{tab:bp and bm operators}. The software implementation of the operations in Table \ref{tab:bp
and bm operators} in \cite{raj:rmtool} uses the companion matrix formulation when $\luvs{1} = \luvs{2}$ and the resultant formulation otherwise. \\
Thus far we have established our ability to encode algebraic distribution as solutions of bivariate polynomial equations and to manipulate the solutions. This sets the stage for defining the class of ``algebraic'' random matrices next.

\section{Class of algebraic random matrices}\label{sec:algebraic random matrices}
 We are interested in identifying canonical random matrix operations for which the limiting eigenvalue distribution of the resulting matrix is an algebraic distribution. This is equivalent to identifying operations for which the transformations in the random matrices can be mapped into transformations of the bivariate polynomial that encodes the limiting eigenvalue distribution function. This motivates the construction of the class of ``algebraic'' random matrices which we shall define next.

The practical utility of this definition, which will become apparent in Section \ref{sec:pol method computational aspects} and \ref{sec:pol method examples} can be succinctly summarized: if a random matrix is shown to be algebraic then its limiting eigenvalue density function can be computed using a simple  root-finding algorithm. Furthermore, if the moments exist, they will satisfy a finite depth linear recursion (see Theorem \ref{th:finite depth recursion}) with polynomial coefficients so that we will often be able to enumerate them efficiently in closed form.  Algebraicity of a random matrix thus acts as a certificate of the computability of its limiting eigenvalue density function and the associated moments. In this chapter our objective is to specify the class of algebraic random matrices by its generators.

\subsection{Preliminaries}
Let ${\bf A}_{N}$, for $N = 1 ,2 , \ldots$ be a sequence of $N \times N$ random matrices with real eigenvalues. Let $F^{{\bf A}_{N}}$ denote the e.d.f., as in (\ref{eq:intro edf}). Suppose $F^{{\bf A}_{N}}(x)$ converges almost surely (or in probability), for every $x$,  to $F^{A}(x)$ as $N \to \infty$, then we say that ${\bf A}_{N} \mapsto A$. We denote the associated (non-random) limiting probability density function by $f_{A}(x)$. 
\begin{notation}[Mode of convergence of the empirical distribution function]
When necessary we highlight the mode of convergence of the underlying distribution function thus: if ${\bf A}_{N} \overset{a.s.}{\longmapsto} A$ then it is shorthand for the statement that the empirical distribution function of ${\bf A}_{N}$ converges almost surely to the distribution function $F^{A}$; likewise ${\bf A}_{N} \overset{p}{\longmapsto} A$ is shorthand for the statement that the empirical distribution function of ${\bf A}_{N}$ converges in probability to the distribution function $F^{A}$.  When the distinction is not made then almost sure convergence is assumed. 
\end{notation}
\begin{remark}
The element $A$ above is not to be interpreted as a matrix. There is no convergence in the sense of an $\infty \times \infty$ matrix. The notation ${\bf A}_{N} \overset{a.s}{\longmapsto} A$ is shorthand for describing the convergence of the associated distribution functions and not of the matrix itself. We think of $A$ as being an (abstract) element of a probability space with distribution function $F^{A}$ and associated density function $f_{A}$.
\end{remark}
\begin{definition}[Atomic random matrix]
If $f_{A} \in \Patom$ then we say that ${\bf A}_{N}$ is an atomic random matrix. We represent this as ${\bf A}_{N} \mapsto A \in \CMatom$ where $\CMatom$ denotes the class of atomic random matrices.
\end{definition}
\begin{definition}[Algebraic random matrix]
 If $f_{A} \in \Palg$ then we say that ${\bf A}_{N}$ is an algebraically characterizable random matrix (often suppressing the word characterizable for brevity). We represent this as ${\bf A}_{N} \longmapsto A \in \CM$ where $\CM$ denotes the class of algebraic random matrices.  Note that, by definition, $\CMatom \subset \CM$.
\end{definition}

\subsection{Key idea used in proving algebraicity preserving nature of a random matrix transformation}
The ability to describe the class of algebraic random matrices and the technique needed to compute the associated bivariate polynomial is at the crux our investigation. In the theorems that follow, we accomplish the former by cataloguing random matrix operations that preserve algebraicity of the limiting distribution.

Our proofs shall rely on exploiting the fact that some random matrix transformations, say ${\bf A}_{N} \longmapsto {\bf B}_{N}$, can be most naturally expressed as transformations of $\lmzs{A} \longmapsto \lmzs{B}$; others as $\lrgs{A} \longmapsto \lrgs{B}$ while some as $\lsys{A} \longmapsto \lsys{B}$. Hence, we manipulate the bivariate polynomials (using the transformations depicted in Figure \ref{fig:all transforms}) to the form needed to apply the appropriate operational law, which we derive as part of the proof, and then reverse the transformations to obtain the bivariate polynomial $\lmzs{B}$. Once we have derived the operational law for computing $\lmzs{B}$ from $\lmzs{A}$, we have established the algebraicity of the limiting eigenvalue distribution of ${\bf B}_{N}$ and we are done. Readers interested in the operational law may skip directly to Section \ref{sec:pol method computational aspects}.

The following property of the convergence of distributions will be invaluable in the proofs that follow .\\

\begin{proposition}[Continuous mapping theorem]\label{th:cmt}
Let ${\bf A}_{N} \longmapsto A$.  Let $f_{A}$ and $\rocat{A}$ denote the corresponding limiting density function and the atomic component of the support, respectively. Consider the mapping $y=h(x)$ continuous everywhere on the real line except on the set of its discontinuities denoted by ${\cal D}_{h}$. If $D_{h} \cap \rocat{A} =  \emptyset$ then ${\bf B}_{N}=h({\bf A}_{N}) \longmapsto B$. The associated non-random distribution function, $F^{B}$ is given by $F^{B}(y)=F^{A}\left(h^{\langle -1 \rangle}(y)\right)$. The associated probability density function is its distributional derivative. 
\end{proposition}
\begin{proof}
This is a restatement of continuous mapping theorem which follows from well-known facts about the convergence of distributions \cite{billingsley99a}.
\end{proof}

\subsection{Deterministic operations}\label{sec:deterministic algebraic}
We first consider some simple deterministic transformations on an algebraic random matrix ${\bf A}_{N}$ that produce an algebraic random matrix ${\bf B}_{N}$. 

\theorembox{
\begin{theorem}\label{th:simple deterministic}
Let ${\bf A}_{N} \mapsto A \in \CM$ and $p$, $q$, $r$, and $s$ be real-valued scalars. Then, 
\[
{\bf B}_{N} = (p\,{\bf A}_{N}+q\,{\bf I}_{N})/(r\,{\bf A}_{N}+s\,{\bf I}_{N}) \mapsto B  \in \CM,
\]
provided $f_{A}$ does not contain an atom at ${-s/r}$  and $r,s$ are not zero simultaneously.
\end{theorem}}
\begin{proof}
Here we have $h(x) = (p\,x+r)/(q\,x+s)$ which is continuous everywhere except at $x =-s/r$ for $s$ and $r$ not simultaneously zero. From Proposition \ref{th:cmt}, unless $f_{A}(x)$ has an atomic component at $-s/r$, ${\bf B}_{N} \mapsto B$. The Stieltjes transform of $F^{B}$ can be expressed as
\begin{equation}\label{eq:BmobiusA mz}
m_{B}(z)=E_{y}\bigg[ \frac{1}{y-z} \bigg]=E_{x}\bigg[ \frac{r\,x+s}{p\,x+q-z(r\,x+s)} \bigg].
\end{equation}
Equation (\ref{eq:BmobiusA mz}) can be rewritten as
\begin{gather}\label{eq:BmobiusA mz 2}
m_{B}(z)=\int \frac{rx+s}{(p-rz)x+(q-sz)}dF^{A}(x)=\frac{1}{p-rz}\int\frac{ rx+s}{x+\frac{q-sz}{p-rz}}dF^{A}(x).
\end{gather}
With some algebraic manipulations, we can rewrite (\ref{eq:BmobiusA mz 2}) as
\begin{gather}
\begin{split}\label{eq:BmobiusA mz 3}
m_{B}(z)&=\beta_{z} \int\frac{rx+s}{x+\alpha_{z}}dF^{A}(x)=\beta_{z} \Bigg( r \int \frac{x}{x+\alpha_{z}}dF^{A}(x)+s\int \frac{1}{x+\alpha_{z}}dF^{A}(x)\Bigg)\\
&=\beta_{z} \Bigg( r \int dF^{A}(x)-r\,\alpha_{z} \int \frac{1}{x+\alpha_{z}}dF^{A}(x)+s\int \frac{1}{x+\alpha_{z}}dF^{A}(x)\Bigg).
\end{split}
\end{gather}
where  $\beta_{z}=1/(p-r\,z)$ and $\alpha_{z}=(q-s\,z)/(p-r\,z)$. Using the definition of the Stieltjes transform and the identity $\int dF^{A}(x) =1$, we can express $m_{B}(z)$ in (\ref{eq:BmobiusA mz 3}) in terms of $m_{A}(z)$ as
\begin{equation}\label{eq:mobius m inter}
  m_{B}(z)=\beta_{z}\,r+(\beta_{z}\,s-\beta\,r\,\alpha_{z})\,m_{A}(-\alpha_{z}).
\end{equation}
Equation (\ref{eq:mobius m inter}) can, equivalently, be rewritten as
\begin{equation}\label{eq:mobius m}
m_{A}(-\alpha_{z})=\dfrac{m_{B}(z)-\beta_{z}\,r}{\beta_{z}\,s-\beta_{z}\,r\,\alpha_{z}}.
\end{equation}
Equation (\ref{eq:mobius m}) can be expressed as an operational law on $\lmzs{A}$ as
\begin{empheq}[
box=\setlength{\fboxrule}{1pt}\fbox]{equation}
\lmz{B}=\lmzs{A}((m-\beta_{z}\,r)/(\beta_{z}\,s-\beta_{z}\,r\,\alpha_{z}),-\alpha_{z}).
\label{eq:mobius pol transform}
\end{empheq}
Since $\lmzs{A}$ exists, we can obtain $\lmzs{B}$ by applying the transformation in (\ref{eq:mobius pol transform}), and clearing the denominator to obtain the irreducible bivariate polynomial consistent with Remark \ref{remark:equivalence class}. Since $\lmzs{B}$ exists, this proves that $f_{B} \in \Palg$ and  ${\bf B}_{N} \mapsto B \in \CM$. 
\end{proof}

Appropriate substitutions for the scalars $p$, $q$, $r$ and $s$ in Theorem \ref{th:simple deterministic} leads to the following Corollary.

\theorembox{
\begin{corollary}\label{cor:general det}
Let ${\bf A}_{N} \mapsto A \in \CM$ and let $\alpha$ be a real-valued scalar. Then,
\begin{enumerate}
\item ${\bf B}_{N} = {\bf A}_{N}^{-1} \mapsto B \in \CM$, provided $f_{A}$ does not contain at atom at $0$,
\item ${\bf B}_{N} = \alpha \, {\bf A}_{N} \mapsto B \in \CM$, 
\item ${\bf B}_{N} = {\bf A}_{N} + \alpha\,{\bf I}_{N} \mapsto B \in \CM$. 
\end{enumerate}
\end{corollary}}

\theorembox{
\begin{theorem}\label{th:transpose}
Let ${\bf X}_{n,N}$ be an $n \times N$ matrix. If ${\bf A}_{N} = {\bf X}_{n,N}{\bf X}_{n,N}^{'} \mapsto A \in \CM$  then 
\[
{\bf B}_{N}={\bf X}_{n,N}^{'}{\bf X}_{n,N} \mapsto B \in \CM.
\]
\end{theorem}
}
\begin{proof}
Here ${\bf X}_{n,N}$ is an $n \times N$ matrix, so that ${\bf A}_{n}$ and ${\bf B}_{N}$ are $n \times n$ and $N \times N$ sized matrices respectively. Let $c_{N} = n/N$. When $c_{N}<1$, ${\bf B}_{N}$ will have $N-n$ eigenvalues of magnitude zero while the remaining $n$ eigenvalues will be identically equal to the eigenvalues of ${\bf A}_{n}$. Thus, the e.d.f. of ${\bf B}_{N}$ is related to the e.d.f. of ${\bf A}_{n}$ as
\begin{gather}\label{eq:transpose cless1 1}
\begin{split}
F^{{\bf B}_{N}}(x) & =\frac{N-n}{N}\,\mathbb{I}_{[0,\infty)}+\frac{n}{N}F^{{\bf A}_{n}}(x)\\
&=(1-c_{N})\,\mathbb{I}_{[0,\infty)}+c_{N}\,F^{{\bf A}_{n}}(x).
\end{split}
\end{gather}
where $\mathbb{I}_{[0,\infty)}$ is the indicator function that is equal to $1$ when $x \geq 0$ and is equal to zero otherwise. 

Similarly, when $c_{N}>1$, ${\bf A}_{n}$ will have $n-N$ eigenvalues of magnitude zero while the remaining $N$ eigenvalues will be identically equal to the eigenvalues of ${\bf B}_{N}$. Thus the e.d.f. of ${\bf A}_{n}$ is related to the e.d.f. of ${\bf B}_{N}$ as
\begin{gather}\label{eq:transpose cmore1 1}
\begin{split}
F^{{\bf A}_{n}}(x)&=\frac{n-N}{n}\,\mathbb{I}_{[0,\infty)}+\frac{N}{n}F^{{\bf B}_{N}}(x)\\
&=\Bigg(1-\frac{1}{c_{N}}\Bigg)\,\mathbb{I}_{[0,\infty)}+\frac{1}{c_{N}}\,F^{{\bf B}_{N}}(x).
\end{split}
\end{gather}
Equation (\ref{eq:transpose cmore1 1}) is  (\ref{eq:transpose cless1 1}) rearranged; so we do not need to differentiate between the case when $c_{N}<1$ and $c_{N}>1$.

Thus, as $n,N \to \infty$ with $c_{N} = n/N \to c$, if $F^{{\bf A}_{n}}$ converges to a non-random d.f. $F^{A}$, then $F^{{\bf B}_{N}}$ will also converge to a non-random d.f. $F^{B}$ related to  $F^{A}$ by
\begin{equation}\label{eq:transpose general c}
F^{B}(x)=(1-c)\mathbb{I}_{[0,\infty)}+c\,F^{A}(x).
\end{equation}
From (\ref{eq:transpose general c}), it is evident that the Stieltjes transform of the limiting distribution functions $F^{A}$ and $F^{B}$ are related as \begin{equation}\label{eq:transpose cmore1 step1}
m_{A}(z)=-\Bigg(1-\frac{1}{c}\Bigg)\frac{1}{z}+\frac{1}{c}m_{B}(z).
\end{equation}
Rearranging the terms on either side of (\ref{eq:transpose cmore1 step1}) allows us to express $m_{B}(z)$ in terms of $m_{A}(z)$ as
\begin{equation}\label{eq:transpose cmore1 m}
m_{B}(z)=-\frac{1-c}{z}+c\,m_{A}(z).
\end{equation}
Equation (\ref{eq:transpose cmore1 m}) can be expressed as an operational law on $\lmzs{A}$ as
\begin{empheq}[
box=\setlength{\fboxrule}{1pt}\fbox]{equation}
\lmz{B}=\lmzs{A}\Bigg(-\Bigg(1-\frac{1}{c}\Bigg)\frac{1}{z}+\frac{1}{c}m,z\Bigg).
\label{eq:transpose pol transform}
\end{empheq}
Given $\lmzs{A}$, we can obtain $\lmzs{B}$ by using (\ref{eq:transpose pol transform}). Hence ${\bf B}_{N} \mapsto B \in \CM$.
\end{proof}

\theorembox{
\begin{theorem}
Let ${\bf A}_{N} \mapsto A \in \CM$. Then $${\bf B}_{N} = \left({\bf A}_{N}\right)^{2} \mapsto B \in \CM.$$ 
\end{theorem}
}
\begin{proof}
Here we have $h(x) = x^2$ which is continuous everywhere. From Proposition \ref{th:cmt}, ${\bf B}_{N} \mapsto B$. The Stieltjes transform of $F^{B}$ can be expressed as
\begin{equation}\label{eq:BsquareA mz}
m_{B}(z)=E_{Y}\bigg[ \frac{1}{y-z} \bigg]=E_{X}\bigg[\frac{1}{x^{2}-z} \bigg].
\end{equation}
Equation (\ref{eq:BsquareA mz}) can be rewritten as
\begin{align}\label{eq:BsquareA mz 2}
m_{B}(z)&=\frac{1}{2\sqrt{z}}\int \frac{1}{x-\sqrt{z}} dF^{A}(x) - \frac{1}{2\sqrt{z}}\int \frac{1}{x+\sqrt{z}} dF^{A}(x) \\
&= \frac{1}{2\sqrt{z}}m_{A}(\sqrt{z})-\frac{1}{2\sqrt{z}}m_{A}(-\sqrt{z}).
\end{align}
Equation (\ref{eq:BsquareA mz 2}) leads to the operational law
\begin{empheq}[
box=\setlength{\fboxrule}{1pt}\fbox]{equation}
\lmz{B}=\lmzs{A}(2m\sqrt{z},\sqrt{z}) \bp{m} \lmzs{A}(-2m\sqrt{z},\sqrt{z}).
\label{eq:square transform}
\end{empheq}
Given $\lmzs{A}$, we can obtain $\lmzs{B}$ by using (\ref{eq:square transform}). This proves that ${\bf B}_{N} \mapsto B \in \CM$.
\end{proof}

\theorembox{
\begin{theorem}\label{th:advanced deterministic}
  Let ${\bf A}_{n} \mapsto A \in  \CM$ and ${\bf B}_{N} \mapsto B \in \CM$. Then, $${\bf C}_{M} = {\rm diag}({\bf A}_{n},{\bf B}_{N}) \mapsto C \in \CM,$$ where $M = n +N$ and $n/N \to c >0$ as $n,N \to \infty$.
\end{theorem}
}
\begin{proof}
Let ${\bf C}_{N}$ be an $N \times N$ block diagonal matrix formed from the $n \times n$ matrix ${\bf A}_{n}$ and the $M \times M$ matrix ${\bf B}_{M}$. Let $c_{N} = n/N$. The e.d.f. of ${\bf C}_{N}$ is given by
\begin{equation*}
F^{{\bf C}_{N}} = c_{N}\,F^{{\bf A}_{n}} + (1 - c_{N})\,F^{{\bf B}_{M}}.
\end{equation*}
Let $n,N \to \infty$ and $c_{N}=n/N \to c$. If $F^{{\bf A}_{n}}$ and $F^{{\bf B}_{M}}$ converge in distribution almost surely (or in probability) to non-random d.f.'s $F^{A}$ and $F^{B}$ respectively, then $F^{{\bf C}_{N}}$ will also converge in distribution almost surely (or in probability) to a non-random distribution function $F^{C}$ given by
\begin{equation}\label{eq:edf block C}
F^{C}(x) = c\,F^{A}(x) + (1-c)\,F^{B}(x).
\end{equation}
The Stieltjes transform of the distribution function $F^{C}$ can hence be written in terms of the Stieltjes transforms of the distribution functions $F^{A}$ and $F^{B}$ as
\begin{equation}\label{eq:A1diagA2 2}
m_{C}(z)=c\,m_{A}(z)+(1-c)\,m_{B}(z)
\end{equation}
Equation (\ref{eq:A1diagA2 2}) can be expressed as an operational law on the bivariate polynomial $\lmz{A}$ as
\begin{empheq}[
box=\setlength{\fboxrule}{1pt}\fbox]{equation} \lmzs{C} =\lmzs{A}\left(\frac{m}{c},z\right) \bp{m}
\lmzs{B}\left(\frac{m}{1-c},z\right).
\end{empheq}
Given $\lmzs{A}$ and $\lmzs{B}$, and the definition of the $\bp{m}$ operator in Section \ref{sec:kronecker}, $\lmzs{C}$ is a polynomial which can be constructed explicitly. This proves that ${\bf C}_{N} \mapsto C \in \CM$.
\end{proof}

\theorembox{
\begin{theorem}\label{th:subspace projection}
If ${\bf A}_{n} = {\rm diag}({\bf B}_{N},\alpha\,{\bf I}_{n-N})$ and $\alpha$ is a real valued scalar. Then,
 $${\bf B}_{N} \mapsto B \in \CM,$$
as $n,N \to \infty$ with $c_{N}=n/N \to c$, 
\end{theorem}
}
\begin{proof}
Assume that as $n,N \to \infty$, $c_{N}=n/N \to c$. As we did in the proof of Theorem \ref{th:advanced deterministic}, we can show that the Stieltjes transform $m_{A}(z)$ can be expressed in terms of $m_{B}(z)$ as
\begin{equation}\label{eq:subspace manipulation m}
m_{A}(z)=\Bigg(\frac{1}{c}-1\Bigg)\frac{1}{\alpha-z}+\frac{1}{c}\,m_{B}(z).
\end{equation}
This allows us to express $\lmz{B}$ in terms of $\lmz{A}$ using the relationship in (\ref{eq:subspace manipulation m}) as
\begin{empheq}[
box=\setlength{\fboxrule}{1pt}\fbox]{equation}
\lmz{B}=\lmzs{A}\Bigg(-\Bigg(\frac{1}{c}-1\Bigg)\frac{1}{\alpha-z}+\frac{1}{c}\,m,z\Bigg).
\label{eq:subspace pol transform}
\end{empheq}
We can hence obtain $\lmzs{B}$ from $\lmzs{A}$ using (\ref{eq:subspace pol transform}). This proves that  ${\bf B}_{N} \mapsto B \in \CM$.
\end{proof}

\theorembox{
\begin{corollary}
Let ${\bf A}_{N} \mapsto A \in \CM$. Then $${\bf B}_{N} ={\rm diag}({\bf A}_{n},\alpha\, {\bf I}_{N-n}) \mapsto B \in \CM,$$ 
for $n/N \to c>0$ as $n,N \to \infty$.
\end{corollary}}
\begin{proof}
This follows directly from Theorem \ref{th:advanced deterministic}.
\end{proof}

\subsection{Gaussian-like operations}
\noindent We now consider some simple stochastic transformations that ``blur'' the eigenvalues of ${\bf A}_{N}$ by injecting additional randomness. We show that canonical operations involving an algebraic random matrix ${\bf A}_{N}$ and Gaussian-like and Wishart-like random matrices (defined next) produce an algebraic random matrix ${\bf B}_{N}$. \\

\begin{definition}[Gaussian-like random matrix]\label{def:gaussianlike matrix}
Let ${\bf Y}_{N,L}$ be an $N \times L$ matrix with independent, identically distributed (i.i.d.) elements having zero mean, unit variance and bounded higher order moments. We label the matrix ${\bf G}_{N,L}=\frac{1}{\sqrt{L}}{\bf Y}_{N,L}$  as a Gaussian-like random matrix.\\
\end{definition}
\noindent We can sample a Gaussian-like random matrix in \matlab as 
\begin{verbatim}
G = sign(randn(N,L))/sqrt(L);
\end{verbatim}
Gaussian-like matrices are labelled thus because they exhibit the same limiting behavior in the $N \to \infty$ limit as ``pure'' Gaussian matrices which may be sampled in \matlab  as 
\begin{verbatim}
G = randn(N,L)/sqrt(L);
\end{verbatim}
\begin{definition}[Wishart-like random matrix]
Let ${\bf G}_{N,L}$ be a Gaussian-like random matrix. We label the matrix $W_{N} = {\bf G}_{N,L} \times {\bf G}_{N,L}'$ as a Wishart-like random matrix. Let $c_{N}=N/L$. We denote a Wishart-like random matrix thus formed by ${\bf W}_{N}(c_{N})$. \\
\end{definition}

\begin{remark}[Algebraicity of Wishart-like random matrices]
The limiting eigenvalue distribution of the Wishart-like random matrix has the Mar\v{c}enko-Pastur density which is an algebraic density since $\lmzs{W}$ exists (see Table \ref{tab:Wishart bivariate}). 
\end{remark}

\begin{proposition}Assume that ${\bf G}_{N,L}$ is an $N \times L$ Gaussian-like random matrix. Let ${\bf A}_{N} \toas A$ be an $N \times N$ symmetric/Hermitian random matrix and ${\bf T}_{L} \toas T$ be an $L \times L$ diagonal atomic random matrix respectively. If ${\bf G}_{N,L}$, ${\bf A}_{N}$ and ${\bf T}_{L}$ are independent then ${\bf B}_{N} = {\bf A}_{N} + {\bf G}_{N,L}^{'}{\bf T}_{L} {\bf G}_{N,L} \toas B$, as $c_{L} = N/L \to c$ for $N, L \to \infty$,. The Stieltjes transform $m_{B}(z)$ of the unique distribution function $F^{B}$ is satisfies the equation
\label{th:mandp}    
\begin{empheq}[
box=\setlength{\fboxrule}{1pt} \fbox]{equation} \label{eq:mandp theorem} m_{B}(z)=m_{A}\bigg(z-c\int\frac{x\,
dF^{T}(x)}{1+x\, m_{B}(z)}\bigg).
\end{empheq}
\end{proposition}
\begin{proof}
This result may be found in Mar\v{c}enko-Pastur\cite{marcenko67a} and Silverstein \cite{silverstein95b}.
\end{proof}\\

\noindent We can reformulate Proposition \ref{th:mandp} to obtain the following result on algebraic random matrices. \\

\theorembox{
\begin{theorem}\label{th:simple stochastic}
Let ${\bf A}_{N}$, ${\bf G}_{N,L}$ and ${\bf T}_{L}$ be defined as in Proposition \ref{th:mandp}. Then $${\bf B}_{N} = {\bf A}_{N}+{\bf G}_{L,N}^{'}{\bf T}_{L}{\bf G}_{L,N} \toas B \in \CM,$$ as $c_{L} = N/L\to c$ for $N, L \to \infty$.
\end{theorem}
}
\begin{proof}
Let ${\bf T}_{L}$ be an atomic matrix with $d$ atomic masses of weight $p_{i}$ and magnitude $\lambda_{i}$ for $i=1,2,\ldots,d$. From Proposition \ref{th:mandp}, $m_{B}(z)$ can be written in terms of $m_{A}(z)$ as
\begin{equation}\label{eq:mandp operator atomic}
m_{B}(z)=m_{A}\Bigg(z-c \sum_{i=1}^{d} \frac{p_{i}\,\lambda_{i}}{1+\lambda_{i} \,m_{B}(z)}\Bigg).
\end{equation}
where we have substituted $F^{T}(x)=\sum_{i=1}^{d}p_{i}\,\mathbb{I}_{[\lambda_{i},\infty)}$ into (\ref{eq:mandp theorem}) with $\sum_{i}p_{i}=1$.

Equation (\ref{eq:mandp operator atomic}) can be expressed as an operational law on the bivariate polynomial $\lmzs{A}$ as
\begin{empheq}[
box=\setlength{\fboxrule}{1pt}\fbox]{equation}
\lmz{B}=\lmzs{A}(m,z-\alpha_{m}).
\label{eq:mandp pol transform}
\end{empheq}
where $\alpha_{m}=c \sum_{i=1}^{d} p_{i}\,\lambda_{i}/(1+\lambda_{i}\,m)$. This proves that ${\bf B}_{N} \toas B \in \CM$. 
\end{proof}\\

\begin{proposition}Assume that ${\bf W}_{N}(c_{N})$ is an $N \times N$ Wishart-like random matrix. Let ${\bf A}_{N} \toas A$ be an $N \times N$ random Hermitian non-negative definite matrix. If ${\bf W}_{N}(c_{N})$ and ${\bf A}_{N}$ are independent, then ${\bf B}_{N} = {\bf A}_{N} \times {\bf W}_{N}(c_{N}) \toas B$ as $c_{N} \to c$. The Stieltjes transform $m_{B}(z)$ of the unique distribution function $F^{B}$ satisfies
\label{th:scm} 
\begin{equation}\label{eq:scm canonical equation}
m_{B}(z)=\int\frac{dF^{A}(x)}{\{1-c-c\,z\,m_{B}(z)\}x-z}.
\end{equation}
\end{proposition}
\begin{proof}
This result may be found in Bai and Silverstein  \cite{bai95a,silverstein95b}.
\end{proof}\\

\noindent We can reformulate Proposition \ref{th:scm} to obtain the following result on algebraic random matrices. \\

\theorembox{
\begin{theorem}
Let ${\bf A}_{N}$ and ${\bf W}_{N}(c_{N})$ satisfy the hypothesis of Proposition \ref{th:scm}. Then, $${\bf B}_{N} = {\bf A}_{N} \times {\bf W}_{N}(c_{N}) \toas B \in \CM,$$ as $c_{N} \to c$.
\end{theorem}
}
\begin{proof}
By rearranging the terms in the numerator and denominator, (\ref{eq:scm canonical equation}) can be rewritten as
\begin{equation}\label{eq:product wishart 2}
m_{B}(z)=\frac{1}{1-c-c\,z\,m_{B}(z)}\int\dfrac{dF^{A}(x)}{x-\frac{z}{1-c-c\,z\,m_{B}(z)}}.
\end{equation}
Let $\alpha_{m,z}=1-c-c\,z\,m_{B}(z)$ so that (\ref{eq:product wishart 2}) can be rewritten as
\begin{equation}\label{eq:product wishart 3}
m_{B}(z)=\frac{1}{\alpha_{m,z}}\int\frac{dF^{A}(x)}{x-(z/\alpha_{m,z})}.
\end{equation}
We can express $m_{B}(z)$ in (\ref{eq:product wishart 3}) in terms of $m_{A}(z)$ as
\begin{equation}\label{eq:product wishart m inter}
m_{B}(z)=\frac{1}{\alpha_{m,z}}\,m_{A}(z/\alpha_{m,z}).
\end{equation}
Equation (\ref{eq:product wishart m inter}) can be rewritten as
\begin{equation}\label{eq:product wishart m}
m_{A}(z/\alpha_{m,z})=\alpha_{m,z}\,m_{B}(z).
\end{equation}
Equation (\ref{eq:product wishart m}) can be expressed as an operational law on the bivariate polynomial $\lmzs{A}$ as
\begin{empheq}[
box=\setlength{\fboxrule}{1pt}\fbox]{equation}
\lmz{B}=\lmzs{A}(\alpha_{m,z}\,m,z/\alpha_{m,z}).
\label{eq:product pol transform}
\end{empheq}
This proves that ${\bf B}_{N} \toas B \in \CM$.
\end{proof}

\begin{proposition}\label{th:silvgirko}
Assume that ${\bf G}_{N,L}$ is an $N \times L$ Gaussian-like random matrix. Let ${\bf A}_{N} \toas A$ be an $N \times N$ symmetric/Hermitian random matrix independent of ${\bf G}_{N,L}$, ${\bf A}_{N}$. Let ${\bf A}_{N}^{1/2}$ denote an $N \times L$ matrix.  If $s$ is a positive real-valued scalar then ${\bf B}_{N} = ({\bf A}_{N}^{1/2} + \sqrt{s}\,{\bf G}_{N,L})({\bf A}_{N}^{1/2} + \sqrt{s}\,{\bf G}_{N,L})^{'} \toas B$, as $c_{L} = N/L \to c$ for $N, L \to \infty$. The Stieltjes transform, $m_{B}(z)$ of the unique distribution function  $F^{B}$ satisfies the equation
\begin{equation} \label{eq:array canonical equation} m_{B}(z)=-\int
\frac{dF^{A}(x)}{z\,\{1+s\,c\,m_{B}(z)\}-\frac{x}{1+s\,c\,m_{B}(z)}+s\,(c-1)}.
\end{equation}
\end{proposition}
\begin{proof}
This result is found in Dozier and Silverstein\cite{silverstein04a}.
\end{proof}

\noindent We can reformulate Proposition \ref{th:silvgirko} to obtain the following result on algebraic random matrices. 

\theorembox{
\begin{theorem}
Assume ${\bf A}_{N}$, ${\bf G}_{N,L}$ and $s$ satisfy the hypothesis of Proposition \ref{th:silvgirko}. Then 
 $${\bf B}_{N} =({\bf A}_{N}^{1/2} + \sqrt{s}\,{\bf G}_{N,L})({\bf A}_{N}^{1/2}+\sqrt{s}\,{\bf G}_{N,L})^{'} \toas B  \in \CM,$$
as $c_{L} = N/L \to c$ for $N,L \to \infty$.
\end{theorem}
}
\begin{proof}
 By rearranging the terms in the numerator and denominator, (\ref{eq:array canonical
equation}) can be rewritten as
\begin{equation}\label{eq:grammian 2}
m_{B}(z)=\int \frac{\{1+s\,c\,m_{B}(z)\}\,dF^{A}(x)}{x-\{1+s\,c\,m_{B}(z)\}(z\,\{1+s\,c\,m_{B}(z)\}+(c-1)\,s)}.
\end{equation}
Let $\alpha_{m}=1+s\,c\,m_{B}(z)$ and $\beta_{m}=\{1+s\,c\,m_{B}(z)\}(z\,\{1+s\,c\,m_{B}(z)\}+(c-1)\,s)$, so that $\beta=\alpha_{m}^{2}\,z+\alpha_{m}\,s(c-1)$. Equation (\ref{eq:grammian 2}) can hence be rewritten as
\begin{equation}\label{eq:grammian 3}
m_{B}(z)=\alpha_{m} \int \frac{dF^{A}(x)}{x-\beta_{m}}.
\end{equation}
Using the definition of the Stieltjes transform in (\ref{eq:mz}), we can express $m_{B}(z)$ in (\ref{eq:grammian 3}) in terms of $m_{A}(z)$ as \begin{gather}\label{eq:grammian m inter}
\begin{split}
m_{B}(z)&=\alpha_{m}\,m_{A}(\beta_{m})\\
&=\alpha_{m}\,m_{A}(\alpha_{m}^{2}\,z+\alpha_{m}(c-1)s).
\end{split}
\end{gather}
Equation (\ref{eq:grammian m inter}) can, equivalently, be rewritten as
\begin{equation}\label{eq:grammian m}
m_{A}(\alpha_{m}^{2}\,z+\alpha_{m}(c-1)s)=\frac{1}{\alpha_{m}}m_{B}(z).
\end{equation}
Equation (\ref{eq:grammian m}) can be expressed as an operational law on the bivariate polynomial $\lmzs{}$ as

\begin{empheq}[
box=\setlength{\fboxrule}{1pt}\fbox]{equation}
\lmz{B}=\lmzs{A}(m/\alpha_{m},\alpha^{2}\,z+\alpha_{m}\,s(c-1)).
\label{eq:array pol transform}
\end{empheq}
This proves that ${\bf B}_{N} \toas B \in \CM$.
\end{proof}

\subsection{Sums and products}\label{sec:sum products of algebraic}

\begin{proposition}\label{th:free add}
Let ${\bf A}_{N} \toinp A$ and ${\bf B}_{N} \toinp B$ be $N \times N$ symmetric/Hermitian random matrices. Let ${\bf Q}_{N}$ be a Haar distributed unitary/orthogonal matrix independent of ${\bf A}_{N}$ and ${\bf B}_{N}$. Then ${\bf C}_{N}={\bf A}_{N} + {\bf Q}_{N}{\bf B}_{N}{\bf Q}_{N}' \toinp C$. The associated distribution function $F^{C}$ is the unique distribution function whose R transform satisfies
\begin{equation}\label{eq:free add}
r_{C}(g) = r_{A}(g)+r_{B}(g).
\end{equation}
\end{proposition}
\begin{proof}
This result was obtained by Voiculescu in \cite{voiculescu86a}.
\end{proof}\\

\noindent We can reformulate Proposition \ref{th:free add} to obtain the following result on algebraic random matrices. 

\theorembox{
\begin{theorem}\label{th:free rtransform rule}
Assume that ${\bf A}_{N}$, ${\bf B}_{N}$ and ${\bf Q}_{N}$ satisfy the hypothesis of Proposition \ref{th:free add}. Then,
$${\bf C}_{N} = {\bf A}_{N} + {\bf Q}_{N}{\bf B}_{N}{\bf Q}_{N}^{'} \toinp C \in \CM$$
\end{theorem}
}
\begin{proof}
Equation (\ref{eq:free add}) can be expressed as an operational law on the bivariate polynomials $\lrgs{A}$ and $\lrgs{B}$ as
\begin{empheq}[
box=\setlength{\fboxrule}{1pt}\fbox]{equation} \lrgs{C} =\lrgs{A} \bp{r} \lrgs{B}
\end{empheq}
If $\lmzs{}$ exists then so does $\lrgs{}$ and vice-versa. This proves that ${\bf C}_{N} \toinp C \in \CM$.
\end{proof}

\begin{proposition}\label{th:free multiply}
Let ${\bf A}_{N} \toinp A$ and ${\bf B}_{N} \toinp B$ be $N \times N$ symmetric/Hermitian random matrices. Let ${\bf Q}_{N}$ be a Haar distributed unitary/orthogonal matrix independent of ${\bf A}_{N}$ and ${\bf B}_{N}$. Then ${\bf C}_{N}={\bf A}_{N} \times {\bf Q}_{N}{\bf B}_{N}{\bf Q}_{N}^{'} \toinp C$ where ${\bf C}_{N}$ is defined only if ${\bf C}_{N}$ has real eigenvalues for every sequence ${\bf A}_{N}$ and ${\bf B}_{N}$. The associated distribution function $F^{C}$ is the unique distribution function whose S transform satisfies
\begin{equation}\label{eq:free multiply}
s_{C}(y) = s_{A}(y)s_{B}(y).
\end{equation}
\end{proposition}
\begin{proof}
This result was obtained by  Voiculescu in \cite{voiculescu87a,voiculescu91a}.
\end{proof}\\

\noindent We can reformulate Proposition \ref{th:free multiply} to obtain the following result on algebraic random matrices. 

\theorembox{
\begin{theorem}\label{th:free stransform rule}
Assume that ${\bf A}_{N}$, and ${\bf B}_{N}$ satisfy the hypothesis of Proposition \ref{th:free multiply}. Then $${\bf C}_{N} = {\bf A}_{N} \times {\bf Q}_{N}{\bf B}_{N}{\bf Q}_{N}^{'} \toinp C \in \CM.$$ 
\end{theorem}
}
\begin{proof}
Equation (\ref{eq:free multiply}) can be expressed as an operational law on the bivariate polynomials $\lsys{A}$ and $\lsys{B}$ as
\begin{empheq}[
box=\setlength{\fboxrule}{1pt}\fbox]{equation} \lsys{C} =\lsys{A} \bt{s} \lsys{B}
\end{empheq}
If $\lmzs{}$ exists then so does $\lsys{}$ and vice versa. This proves that ${\bf B}_{N} \toinp B \in \CM$.
\end{proof}

\begin{definition}[Orthogonally/Unitarily invariant random matrix]
If the joint distribution of the elements of a random matrix ${\bf A}_{N}$ is invariant under orthogonal/unitary transformations, it is referred to as an orthogonally/unitarily invariant random matrix.\\
\end{definition}
\noindent If ${\bf A}_{N}$ (or ${\bf B}_{N}$) or both are an orthogonally/unitarily invariant sequences of random matrices then Theorems \ref{th:free rtransform rule} and \ref{th:free stransform rule} can be stated more simply. \\

\theorembox{
\begin{corollary}\label{th:free convolution invariant}
Let ${\bf A}_{N} \toinp A \in \CM$ and ${\bf B}_{N} \to B \toinp \CM$ be a orthogonally/unitarily invariant random matrix independent of ${\bf A}_{N}$. Then,
\begin{enumerate}
\item ${\bf C}_{N} = {\bf A}_{N} + {\bf B}_{N} \toinp C \in \CM$
\item ${\bf C}_{N} = {\bf A}_{N} \times {\bf B}_{N} \toinp C \in \CM$  
\end{enumerate}
Here multiplication is defined only if ${\bf C}_{N}$ has real eigenvalues for every sequence ${\bf A}_{N}$ and ${\bf B}_{N}$.\\
\end{corollary}
}
\noindent When both the limiting eigenvalue distributions of ${\bf A}_{N}$ and ${\bf B}_{N}$ have compact support, it is possible to strengthen the mode of convergence in Theorems \ref{th:free rtransform rule} and \ref{th:free stransform rule} to almost surely \cite{hiai00a}. We suspect that almost sure convergence must hold when the distributions are not compactly supported; this remains an open problem.

%\Chapter{The polynomial method: Computational aspects}\label{chap:pol method computational aspects}

\section{Operational laws on bivariate polynomials}\label{sec:pol method computational aspects}
The key idea behind the definition of algebraic random matrices in Section \ref{sec:algebraic random matrices} was that when the limiting eigenvalue distribution of a random matrix can be encoded by a bivariate polynomial, then for the broad class of random matrix operations identified in Section \ref{sec:algebraic random matrices}, algebraicity of the eigenvalue distribution is preserved under the transformation.

These operational laws, the associated random matrix transformation and the symbolic \matlab code for the operational law are summarized in Tables \ref{tab:summary bivariate transforms}-\ref{tab:remaining transformations}. The remainder of this chapter discusses techniques for extracting the density function from the polynomial and the special structure in the moments that allows them to be efficiently enumerated using symbolic methods.

\begin{table}[p]
\centering
\captionstyle{flushleft}
\includegraphics[width=6.35in,angle=90]{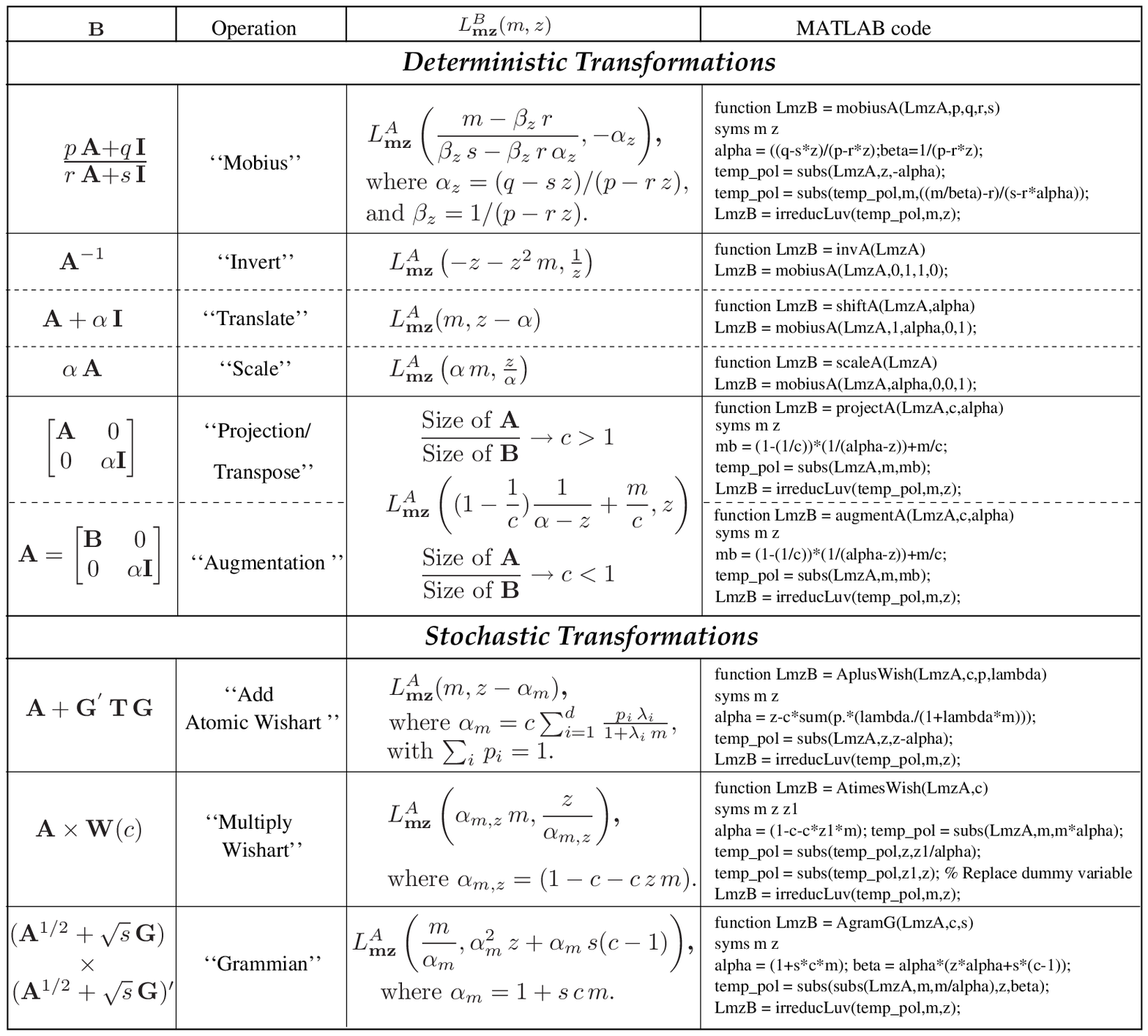}
\caption{Operational laws on the bivariate polynomial encodings (and their computational realization in \matlab) corresponding to a class of deterministic and stochastic transformations. The Gaussian-like random matrix ${\bf G}$ is an $N \times L$, the Wishart-like matrix ${\bf W}(c)={\bf G}\,{\bf G}^{'}$ where $N/L \to c >0$ as $N,L \to \infty$, and the matrix ${\bf T}$ is a diagonal atomic random matrix.}
\label{tab:summary bivariate transforms}
\end{table}\index{polynomial method!operational law}\index{polynomial method!MATLAB code}\index{polynomial method!computational software}

\begin{table}
 \centering
  \subtable[$\lmzs{A} \longmapsto \lmzs{B}$ for  ${\bf A} \longmapsto {\bf B} = {\bf A}^{2}$.]
{ \label{tab:B to Asq}
\begin{tabular}{|c|c|}
\hline
Operational Law & \matlab Code\\
\hline
{\setlength{\tabcolsep}{0pt}
\begin{tabular}{ccccc}
&& $\lmzs{A}$ &&\\
& $\swarrow$ && $\searrow$ &\\
$\lmzs{A}(2m\sqrt{z},\sqrt{z})$&&&&$\lmzs{A}(-2m\sqrt{z},-\sqrt{z})$\\
& $\searrow$ && $\swarrow$ &\\
&& $\,\:\bp{m}$ &&\\
&& $\downarrow$ &&\\
&& $\lmzs{B}$ &&\\
\end{tabular}}&
\begin{small}
\begin{tabular}{l}
\texttt{function LmzB = squareA(LmzA)}\\
\texttt{syms m z }\\
\texttt{}\\
\texttt{Lmz1 = subs(LmzA,z,sqrt(z));}\\
\texttt{Lmz1 = subs(Lmz1,m,2*m*sqrt(z));}\\
\texttt{Lmz2 = subs(LmzA,z,-sqrt(z));}\\
\texttt{Lmz2 = subs(Lmz2,m,-2*m*sqrt(z));}\\
\texttt{}\\
\texttt{LmzB = L1plusL2(Lmz1,Lmz2,m);}\\
\texttt{LmzB = irreducLuv(LmzB,m,z);}\\
\end{tabular}
\end{small}\\
\hline
\end{tabular}
}      \\[0.5in]          
\subtable[$\lmzs{A}, \lmzs{A} \longmapsto \lmzs{C}$ for ${\bf A}, {\bf B} \longmapsto {\bf C} =
\textrm{diag}({\bf A},{\bf B})$ where Size of ${\bf A}$/ Size of ${\bf C} \to c$.]{\label{tab:B to A1diagA2}
\begin{tabular}{|c|c|}\hline
Operational Law & \matlab Code\\ \hline
 {\setlength{\tabcolsep}{0pt}
  \begin{tabular}{ccccc}
    $\lmzs{A}$   &&&& $\lmzs{B}$\\
    $\downarrow$ &&&& $\downarrow$\\
    $\lmzs{A}(\frac{m}{c},z)$ &&&& $\lmzs{B}(\frac{m}{1-c},z)$\\
    & $\searrow$ && $\swarrow$ &\\
    && $\:\,\bp{m}$ &&\\
    && $\downarrow$ &&\\
    && $\lmzs{C}$ &&\\
  \end{tabular}}
&  \begin{small}
    \begin{tabular}{l}
    \texttt{function LmzC = AblockB(LmzA,LmzB,c)}\\
    \texttt{syms m z mu}\\
    \texttt{} \\
    \texttt{LmzA1 = subs(LmzA,m,m/c);}\\
    \texttt{LmzB1 = subs(LmzB,m,m/(1-c));}\\
    \texttt{}\\
    \texttt{LmzC = L1plusL2(LmzA1,LmzB1,m);}\\
    \texttt{LmzC = irreducLuv(LmzC,m,z);}\\
  \end{tabular}
   \end{small}\\
\hline
\end{tabular}
} \\
\caption{Operational laws on the bivariate polynomial encodings for some deterministic random matrix transformations. The operations $\bp{u}$ and $\bt{u}$ are defined in Table \ref{tab:bp and bm operators}.} 
\label{tab:remaining deterministic transformations}
\end{table}\index{polynomial method!operational law}\index{polynomial method!MATLAB code}\index{polynomial method!computational software}

\begin{table}
\centering
\subtable[$\lmzs{A}, \lmzs{B} \longmapsto \lmzs{C}$ for ${\bf A}, {\bf B} \longmapsto {\bf C} =  {\bf A} + {\bf Q} {\bf B} {\bf Q}^{'}$.]{\label{tab:A1plusA2}
\begin{tabular}{|c|c|}\hline
Operational Law & \matlab Code\\ \hline
 {\setlength{\tabcolsep}{0pt}
  \begin{tabular}{ccccc}
    $\lmzs{A}$   &&&& $\lmzs{B}$\\
    $\downarrow$ &&&& $\downarrow$\\
    $\lrgs{A}$   &&&& $\lrgs{B}$\\
    & $\searrow$ && $\swarrow$ &\\
    && $\:\,\bp{r}$ &&\\
    && $\downarrow$ &&\\
    && $\lrgs{C}$ &&\\
    && $\downarrow$ &&\\
    && $\lmzs{C}$ &&\\
  \end{tabular}}
& \begin{small}
  \begin{tabular}{l}
    \texttt{function LmzC = AplusB(LmzA,LmzB)}\\
    \texttt{syms m z r g}\\
    \texttt{}\\
    \texttt{LrgA = Lmz2Lrg(LmzA);}\\
    \texttt{LrgB = Lmz2Lrg(LmzB);}\\
    \texttt{}\\
    \texttt{LrgC = L1plusL2(LrgA,LrgB,r);}\\
    \texttt{}\\
    \texttt{LmzC = Lrg2Lmz(LrgC);}\\
  \end{tabular}
  \end{small}
  \\
\hline
\end{tabular}
} \\[0.5in]
\subtable[$\lmzs{A}, \lmzs{B} \longmapsto \lmzs{C}$ for ${\bf A}, {\bf B} \longmapsto {\bf C} = {\bf A} \times {\bf Q} {\bf B} {\bf Q}^{'}$.]{\label{tab:A1timesA2}
\begin{tabular}{|c|c|}\hline
Operational Law & \matlab Code\\ \hline
 {\setlength{\tabcolsep}{0pt}
  \begin{tabular}{ccccc}
    $\lmzs{A}$   &&&& $\lmzs{B}$\\
    $\downarrow$ &&&& $\downarrow$\\
    $\lsys{A}$   &&&& $\lsys{B}$\\
    & $\searrow$ && $\swarrow$ &\\
    && $\:\,\bt{s}$ &&\\
    && $\downarrow$ &&\\
    && $\lsys{C}$ &&\\
    && $\downarrow$ &&\\
    && $\lmzs{C}$ &&\\
  \end{tabular}}
& \begin{small}
  \begin{tabular}{l}
    \texttt{function LmzC = AtimesB(LmzA,LmzB)}\\
    \texttt{syms m z s y} \\
    \texttt{}\\
    \texttt{LsyA = Lmz2Lsy(LmzA);}\\
    \texttt{LsyB = Lmz2Lsy(LmzB);}\\
    \texttt{}\\
    \texttt{LsyC = L1timesL2(LsyA,LsyB,s);}\\
    \texttt{}\\
    \texttt{LmzC = Lsy2Lmz(LsyC);}\\
  \end{tabular}
  \end{small}
  \\
\hline
\end{tabular}
}
\caption{Operational laws on the bivariate polynomial encodings for some canonical random matrix transformations. The operations $\bp{u}$ and $\bt{u}$ are defined in Table \ref{tab:bp and bm operators}.} 
\label{tab:remaining transformations}
\end{table}\index{polynomial method!operational law}\index{polynomial method!MATLAB code}\index{polynomial method!computational software}\index{polynomial method!free probability}

\newpage

\section{Interpreting the solution curves of polynomial equations}\label{sec:level density}

Consider a bivariate polynomial $\lmzs{}$. Let $\Dm$ be the degree of $\lmz{}$ with respect to $m$ and $l_{k}(z)$, for $k=0,\ldots, \Dm$, be polynomials in $z$ that are the coefficients of $m^{k}$. For every $z$ along the real axis, there are at most $\Dm$ solutions to the polynomial equation $\lmz{}=0$. The solutions of the bivariate polynomial equation $\lmzs{}=0$ define a locus of points $(m,z)$ in $\mathbb{C} \times \mathbb{C}$ referred to as a complex algebraic curve. Since the limiting density is over $\mathbb{R}$, we may focus on real values of $z$.\index{algebraic curves}

For almost every $ z \in \mathbb{R}$, there will be $\Dm$ values of $m$. The exception consists of the singularities of $\lmz{}$. A singularity occurs at $z=z_{0}$ if:
\begin{itemize}
\item There is a reduction in the degree of $m$ at $z_{0}$ so that there are less than $\Dm$ roots for $z=z_{0}$. This occurs when $l_{\Dm}(z_{0})=0$. Poles of $\lmz{}$ occur if some of the $m$-solutions blow up to infinity.
\item There are multiple roots of $\lmzs{}$ at $z_{0}$ so that some of the values of $m$ coalesce. 
\end{itemize}

The singularities constitute the so-called exceptional set of $\lmz{}$. Singularity analysis, in the context of algebraic functions, is a well studied problem \cite{flajolet01a} from which we know that the singularities of $\lmz{A}$ are constrained to be \textit{branch points}. \index{algebraic curves!singularities}

A \textit{branch} of the algebraic curve $\lmz{}=0$ is the choice of a locally analytic function $m_{j}(z)$ defined outside the exceptional set of $\lmz{A}$ together with a connected region of the $\mathbb{C} \times \mathbb{R}$ plane throughout which this particular choice $m_{j}(z)$ is analytic. These properties of singularities and branches of algebraic curve are helpful in determining the atomic and non-atomic component of the encoded probability density from $\lmzs{}$. We note that, as yet, we do not have a fully automated algorithm for extracting the limiting density function from the bivariate polynomial. Development of efficient computational algorithms that exploit the algebraic properties of the solution curve would be of great benefit to the community.
\index{algebraic curves!singularities!eigenvalue density}

%\begin{notation}[Region of support]\label{not:roc}
%Let $F^{A}(x)$ be a distribution function with $f_{A}(x)$ its distributional derivative. Let $\roc{A}$ denote the support of $f_{A}(x)$, with %$\rocat{A}$ and $\rocnon{A}$ denoting the atomic and non-atomic components of the support, respectively. Thus $\roc{A} = \rocat{A} \cup %\rocnon{A}$.\\
%\end{notation}

\subsection{The atomic component}\index{algebraic curves!singularities!eigenvalue density}
\index{eigenvalue density function!atomic component}
If there are any atomic components in the limiting density function, they will necessarily manifest themselves as poles of $\lmz{}$. This follows from the definition of the Stieltjes transform in (\ref{eq:mz}). As mentioned in the discussion on the singularities of algebraic curves, the poles are located at the roots of $l_{\Dm}(z)$. These may be computed in \maple using the sequence of commands:

\begin{small}
\begin{verbatim}
   > Dm := degree(LmzA,m);
   > lDmz := coeff(LmzA,m,Dm);
   > poles := solve(lDmz=0,z);
\end{verbatim}
\end{small}

We can then compute the Puiseux expansion about each of the poles at $z=z_{0}$. This can be computed in \maple using the \texttt{algcurves} package as:

\begin{small}
\begin{verbatim}
   > with(algcurves):
   > puiseux(Lmz,z=pole,m,1);
\end{verbatim}
\end{small}

For the pole at $z=z_{0}$, we inspect the Puiseux expansions for branches with leading term $1/(z_{0}-z)$. An atomic component in the limiting spectrum occurs if and  only if the coefficient of such a branch is non-negative and not greater than one. This constraint ensures that the branch is associated with the Stieltjes transform of a valid probability distribution function.

Of course, as is often the case with algebraic curves, pathological cases can be easily constructed. For example, more than one branch of the Puiseux expansion might correspond to a candidate atomic component, \ie, the coefficients are non-negative and not greater than one. In our experimentation, whenever this has happened it has been possible to eliminate the spurious branch by matrix theoretic arguments. Demonstrating this rigorously using analytical arguments remains an open problem. 

Sometimes it is possible to encounter a double pole at $z=z_{0}$ corresponding to two admissible weights. In such cases, empirical evidence suggests that the branch with the largest coefficient (less than one) is the ``right'' Puiseux expansion though we have no theoretical justification for this choice.

\subsection{The non-atomic component}\index{algebraic curves!singularities!eigenvalue density}\index{eigenvalue density function!non-atomic component}

The probability density function can be recovered from the Stieltjes transform by applying the inversion formula in (\ref{eq:inversion formula}). Since the Stieltjes transform is encoded in the bivariate polynomial $\lmzs{}$, we accomplish this by first computing all $\Dm$ roots along $z \in \mathbb{R}$ (except at poles or singularities). There will be $\Dm$ roots of which one solution curve will be the ``correct'' solution , \ie, the non-atomic component of the desired density function is the imaginary part of the correct solution normalized by $\pi$.  In \matlab, the $\Dm$ roots can be computed using the sequence of commands:

\begin{small}
\begin{verbatim}
   Lmz_roots = [];
   x_range = [x_start:x_step:x_end];
   for x = x_range
      Lmz_roots_unnorm = roots(sym2poly(subs(Lmz,z,x)));
      Lmz_roots = [Lmz_roots; 
                   real(Lmz_roots_unnorm) + i*imag(Lmz_roots_unnorm)/pi];
   end
\end{verbatim}
\end{small}

The density of the limiting eigenvalue distribution function can be, generically, be expressed in closed form  when $\Dm=2$. When using root-finding algorithms, for $\Dm=2,3$, the correct solution can often be easily identified; the imaginary branch will always appear with its complex conjugate. The density is just the scaled (by $1/\pi$) positive imaginary component. 

When $\Dm \geq 4$, except when $\lmzs{}$ is bi-quadratic for $\Dm = 4$, there is no choice but to manually isolate the correct solution among the numerically computed $\Dm$ roots of the polynomial $\lmzs(m,z)$ at each $z=z_{0}$. The class of algebraic random matrices whose eigenvalue density function can be expressed in closed form is thus a much smaller subset of the class of algebraic random matrices. When the underlying density function is compactly supported, the boundary points will be  singularities of the algebraic curve. 

In particular, when the probability density function is compactly supported and the boundary points are not poles, they occur at points where some values of $m$ coalesce. These points are the roots of the discriminant of $\lmzs{}$, computed in \maple as:

\begin{small}
\begin{verbatim}
   > PossibleBoundaryPoints = solve(discrim(Lmz,m),z);
\end{verbatim}
\end{small}

We suspect that ``nearly all''  algebraic random matrices with compactly supported eigenvalue distribution will exhibit a square root type behavior near boundary points at which there are no poles. In the generic case, this will occur whenever the boundary points correspond to locations where two branches of the algebraic curve coalesce.

For a class of random matrices that includes a subclass of algebraic random matrices, this has been established in \cite{silverstein95a}. This endpoint behavior has also been observed orthogonally/unitarily invariant random matrices whose distribution has the element-wise joint density function of the form

\begin{equation*}
f({\bf A}) = C_{N} \exp\left(-N \Tr \,V({\bf A})\right) d{\bf A} 
\end{equation*}

\noindent where $V$ is an even degree polynomial with positive leading coefficient and $d{\bf A}$ is the Lebesgue measure on $N \times N$ symmetric/Hermitian matrices. In \cite{deift98a}, it is shown that these random matrices have a limiting mean eigenvalue density in the $N \to \infty$ limit that is algebraic and compactly supported. The behavior at the endpoint typically vanishes like a square root, though higher order vanishing at endpoints is possible and a full classification is made in \cite{deift99a}. In \cite{kuijlaars00a} it is shown that square root vanishing is generic. A similar classification for the general class of algebraic random matrices remains an open problem. This problem is of interest because of the intimate connection between the endpoint behavior and the Tracy-Widom distribution. Specifically, we conjecture that ``nearly all'' algebraic random matrices with compactly supported eigenvalue distribution whose density function vanishes as the square root at the endpoints will, with appropriate re-centering and rescaling, exhibit Tracy-Widom fluctuations.

\index{eigenvalue density function!endpoint behavior!conjecture}

Whether the encoded distribution is compactly supported or not, the $-1/z$ behavior of the real part of Stieltjes transform (the principal value) as $z \to \pm \infty$  helps isolate the correct solution. In our experience,  while multiple solution curves might exhibit this behavior, invariably only one solution will have an imaginary branch that, when normalized, will correspond to a valid probability density. Why this always appears to be the case for the operational laws described is a bit of a mystery to us.\\

\begin{figure}[hbp]
 \centering
    \subfigure[Real component. The singularity at zero corresponds to an atom of weight $1/2$. The branch points at $(1\pm \sqrt{2})^{2}$ correspond to the boundary points of the region of support.]{
    \label{fig:mandp real}
    \psfrag{x}{z}
     \includegraphics[scale=0.80]{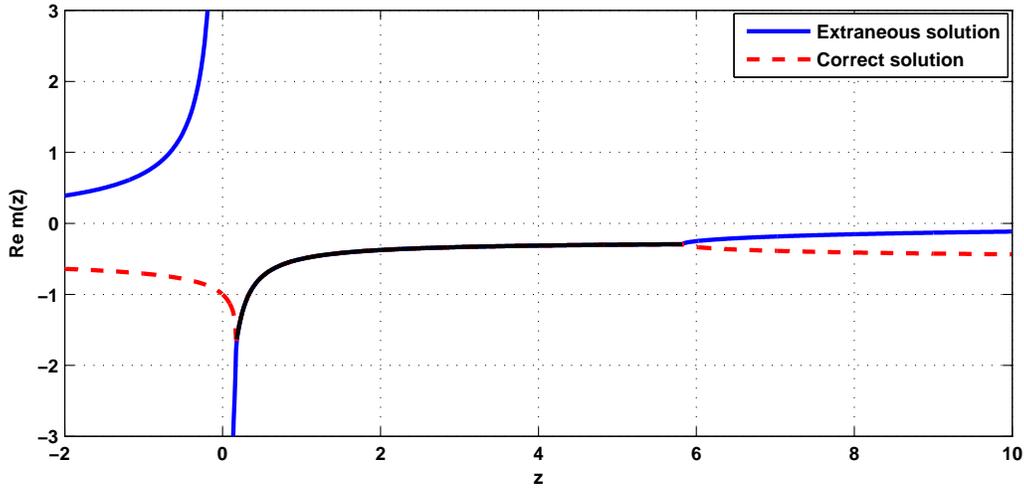}
    
} \\
     \vspace{.3in}
     \subfigure[Imaginary component normalized by $\pi$. The positive component corresponds to the encoded probability density function.]{
     \label{fig:mandp imag}
     \includegraphics[scale=0.80]{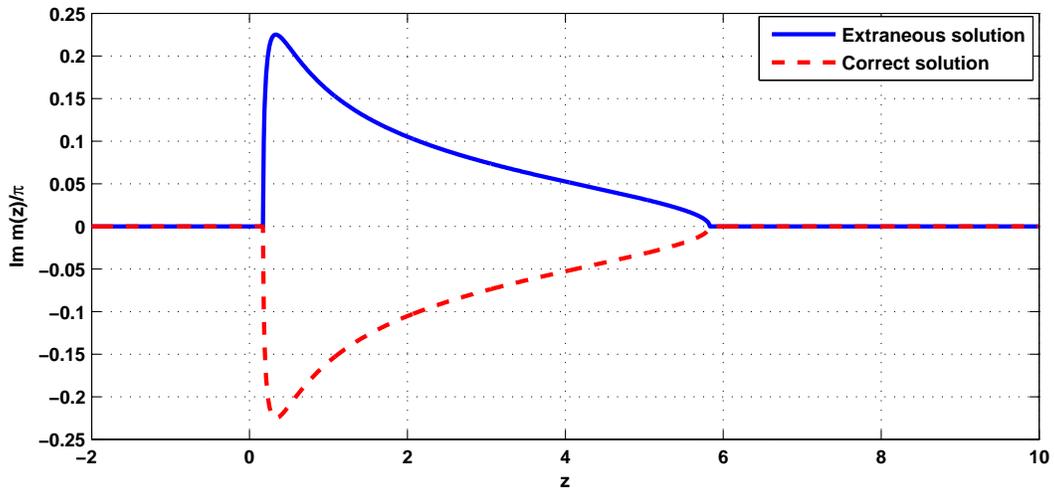}
     } \\
\caption{The real and imaginary components of the algebraic curve defined by the equation $\lmz{}=0$, where $\lmzs{}\equiv cz{m}^{2}- \left( 1-c-z \right) m+1$, which encodes the Mar\v{c}enko-Pastur density. The curve is plotted for $c=2$.  The $-1/z$ behavior of the real part of the ``correct solution'' as $z \to \infty$ is the generic behavior exhibited by the real part of the Stieltjes transform of a valid probability density function. }
\label{fig:algebraic curve example}
\end{figure}

\noindent \textbf{Example}: Consider the Mar\v{c}enko-Pastur density encoded by $\lmzs{}$ given in Table \ref{tab:Wishart bivariate}. The Puiseux expansion about the pole at $z=0$ (the only pole!), has coefficient $(1-1/c)$ which corresponds to an atom only when $c>1$ (as expected using a matrix theoretic argument). Finally, the branch points at $(1\pm \sqrt{c})^{2}$ correspond to boundary points of the compactly supported probability density. Figure \ref{fig:algebraic curve example} plots the real and imaginary parts of the algebraic curve for $c=2$.

\section{Enumerating the moments and free cumulants}\label{sec:moments}\index{moment power series!enumeration}\index{free cumulant power series!enumeration}
In principle, the moments generating function can be extracted from $\lmuzs{}$ by a Puiseux expansion of the algebraic function $\mu(z)$ about $z=0$. When the moments of an algebraic probability distribution exist, there is additional structure in the moments and free cumulants that allows us to enumerate them efficiently. For an algebraic probability distribution, we conjecture that the moments of all order exist if and only if the distribution is compactly supported.\index{moments!existence!conjecture}

\begin{definition}[Rational generating function]\index{rational power series!definition}\label{def:rational gfun}
Let $\mathbb{R}[[x]]$ denote the ring of formal power series (or generating functions) in $x$ with real coefficients. A formal power series (or generating function) $v \in \mathbb{R}[[u]]$ is said to be rational if there exist polynomials in $u$, $P(u)$ and $Q(u)$, $Q(0) \neq 0$ such that 
\[
v(u) = \dfrac{P(u)}{Q(u)}.
\]
\end{definition}

\begin{definition}[Algebraic generating function]\index{algebraic power series!definition}\label{def:algebraic gfun}
Let $\mathbb{R}[[x]]$ denote the ring of formal power series (or generating functions) in $x$ with real coefficients. A formal power series (or generating function) $v \in \mathbb{R}[[u]]$ is said to be algebraic if there exist polynomials in $u$, $P_{0}(u), \ldots, P_{\Du}(u)$, not all identically zero, such that
\[
P_{0}(u)+ P_{1}(u)v+\ldots+P_{\Dv}(u)v^{\Dv}=0.
\]
The degree of $v$ is said to be $\Dv$.
\end{definition}

\begin{definition}[D-finite generating function]\label{def:dfinite}\index{D-finite generating functions!definition}
Let $v \in \mathbb{R}[[u]]$. If there exist polynomials $p_{0}(u), \ldots, p_{d}(u)$, such that
\begin{equation}\label{eq:dfinite}
p_{d}(u)v^{(d)}+p_{d-1}(u)v^{(d-1)}+\ldots+p_{1}(u)v^{(1)}+p_{0}(u)=0,
\end{equation}
where $v^{(j)}=d^{j}v/du^{j}$. Then we say that $v$ is a D-finite (short for differentiably finite) generating function (or power series).  The generating function, $v(u)$, is also referred to as a holonomic function. 
\end{definition}

\begin{definition}[P-recursive coefficients]\index{polynomially recursive coefficients!definition}
Let $a_{n}$ for $n \geq 0$ denote the coefficients of a D-finite series $v$. If there exist polynomials $P_{0},\ldots, P_{e} \in \mathbb{R}[n]$ with $P_{e} \neq 0$, such that
\[
P_{e}(n)a_{n+e}+P_{e-1}(n)a_{n+e-1}+\ldots+P_{0}(n)a_{n}=0,
\]
for all $n \in \mathbb{N}$, then the coefficients $a_{n}$ are said to be P-recursive (short for polynomially recursive). 
\end{definition}

\begin{proposition}\label{th:dfinite theorem}\index{polynomially recursive coefficients!algebraic power series}\index{algebraic power series!polynomially recursive coefficients}\index{D-finite power series!algebraic power series}
Let $v \in \mathbb{R}[[u]]$ be an algebraic power series of degree $\Dv$. Then $v$ is D-finite and satisfies an equation (\ref{eq:dfinite}) of order $\Dv$.
\end{proposition}
\begin{proof}
A proof appears in Stanley \cite[pp.187]{stanley99a}.
\end{proof}

\vspace{0.15in}
\noindent The structure of the limiting moments and free cumulants associated with algebraic densities is described next. 
\vspace{0.15in}

\theorembox{
\begin{theorem}\index{polynomial recursive coefficients!moments}\index{polynomial recursive coefficients!free cumulants}\index{algebraic power series!moment power series}\index{algebraic power series!free cumulant power series}\index{D-finite generating functions!moment power series}\index{D-finite generating functions!free cumulant power series}
\label{th:finite depth recursion}
If $f_{A} \in \Palg$, and the moments exist, then the moment and free cumulant generating functions are algebraic power series. Moreover, both generating functions are D-finite and the coefficients are P-recursive.
\end{theorem}
}
\begin{proof}
If $f_{A} \in \Palg$,  then $\lmzs{A}$ exists. Hence $\lmuzs{A}$ and $\lrgs{A}$ exist, so that $\mu_{A}(z)$ and $r_{A}(g)$ are algebraic power series. By Theorem \ref{th:dfinite theorem} they are D-finite; the moments and free cumulants are hence P-recursive. 
\end{proof}\\

There are powerful symbolic tools available for enumerating the coefficients of algebraic power series. The \maple based package \texttt{gfun} is one such example \cite{salvy94a}. From the bivariate polynomial $\lmuzs{}$, we can obtain the series expansion up to degree \texttt{expansion\_degree} by using the commands:

\begin{small}
\begin{verbatim}
   > with(gfun):
   > MomentSeries = algeqtoseries(Lmyuz,z,myu,expansion_degree,'pos_slopes');
\end{verbatim}
\end{small}

The option \texttt{pos\_slopes} computes only those branches tending to zero. Similarly, the free cumulants can be enumerated from $\lrgs{}$ using the commands:

\begin{small}
\begin{verbatim}
   > with(gfun):
   > FreeCumulantSeries = algeqtoseries(Lrg,g,r,expansion_degree,'pos_slopes');
\end{verbatim}
\end{small}

\noindent For computing expansions to a large order, it is best to work with the recurrence relation. For an algebraic power series $v(u)$, the first \texttt{number\_of\_terms} coefficients can be computed from $\luvs{}$ using the sequence of commands:

\begin{small}
\begin{verbatim}
   > with(gfun):
   > deq := algeqtodiffeq(Luv,v(u));
   > rec := diffeqtorec(deq,v(u),a(n));
   > p_generator := rectoproc(rec,a(n),list):
   > p_generator(number_of_terms);
\end{verbatim}
\end{small} 

\vspace{0.1cm}
\noindent\textbf{Example}: Consider the Mar\v{c}enko-Pastur density encoded by the bivariate polynomials listed in Table \ref{tab:Wishart bivariate}. Using the above sequence of commands, we can enumerate the first five terms of its  moment generating function as 

\begin{small}
\[
\mu(z)=1+z+ \left( c+1 \right) {z}^{2}+ \left( 3\,c+{c}^{2}+1 \right) {z}^{
3}+ \left( 6\,{c}^{2}+{c}^{3}+6\,c+1 \right) {z}^{4}+O \left( {z}^{5}
 \right).
\]
\end{small}

\noindent The moment generating function is a D-Finite power series and satisfies the second order differential equation

\begin{small}
\[
-z+zc-1+ \left( -z-zc+1 \right) \mu \left( z \right) + \left( {z
}^{3}{c}^{2}-2\,{z}^{2}c-2\,{z}^{3}c+z-2\,{z}^{2}+{z}^{3} \right) {
\frac {d}{dz}}\mu \left( z \right)=0,
\]
\end{small}

\noindent with initial condition $\mu(0)=1$. The moments $M_{n} = a(n)$ themselves are P-recursive satisfying the finite depth recursion

\begin{small}
\[
 \left( -2\,c+{c}^{2}+1 \right) na \left( n \right) + \left( 
 \left( -2-2\,c \right) n-3\,c-3 \right) a \left( n+1 \right) +
 \left( 3+n \right) a \left( n+2 \right)=0 
\]
\end{small}

\noindent with the initial conditions $a \left( 0 \right) =1$ and $a\left( 1 \right) =1$. The free cumulants can be analogously computed.

What we find rather remarkable is that for algebraic random matrices, it is often possible to enumerate the moments in closed form even when the limiting density function cannot. The linear recurrence satisfied by the moments may be used to analyze their asymptotic growth.  

When using the sequence of commands described, sometimes more than one solution might emerge. In such cases, we have often found that one can identify the correct solution by checking for the positivity of even moments or the condition $\mu(0)=1$. More sophisticated arguments might be needed for pathological cases. It might involve verifying, using techniques such as those in \cite{akhiezer65a}, that the coefficients enumerated correspond to the moments a valid distribution function.

\section{Computational free probability}\label{sec:computational free probability}

\subsection{Moments of random matrices and asymptotic freeness}\label{sec:moments and freeness}

Assume we know the eigenvalue distribution of two matrices ${\bf A}$ and ${\bf B}$. In general, using that information alone, we cannot say much about the
the eigenvalue distribution of the sum ${\bf A}+{\bf B}$ of the matrices since eigenvalues of the sum of the matrices depend on the eigenvalues of ${\bf A}$ and the eigenvalues of ${\bf B}$, and also on the relation between the eigenspaces of ${\bf A}$ and of ${\bf B}$. However, if we pose this question in the context of $N\times N$-random matrices, then in many situations the answer becomes deterministic in the limit $N\to\infty$. Free probability provides the analytical framework for characterizing this limiting behavior.

\begin{definition}\label{def:limiteigenval}
Let ${\bf A}=({\bf A}_N)_{N\in\NN}$ be a sequence of $N\times N$-random matrices. We say that ${\bf A}$ has a
limit eigenvalue distribution if the limit of all moments
$$\alpha_n:=\lim_{N\to\infty} E[\tr({\bf A}_N^n)]\qquad (n\in\NN)$$
exists, where $E$ denotes the expectation and $\tr$ the normalized trace.
\end{definition}

Using the language of limit eigenvalue distribution as in Definition ~\ref{def:limiteigenval}, our question becomes: Given two random matrix ensembles of $N\times N$-random
matrices, ${\bf A}=({\bf A}_N)_{N\in\NN}$ and ${\bf B}=({\bf B}_N)_{N\in\NN}$, with limit eigenvalue distribution, does their sum ${\bf C}=({\bf C}_N)_{N\in\NN}$, with ${\bf C}_N={\bf A}_N+{\bf B}_N$, have a limit eigenvalue distribution, and furthermore, can we calculate the limit moments $\alpha_n^C$ of ${\bf C}$ out of the limit
moments $(\alpha_k^{\bf A})_{k\geq 1}$ of ${\bf A}$ and the limit moments $(\alpha_k^B)_{k\geq 1}$ of ${\bf B}$ in a deterministic way. It turns out that this is the case if the two ensembles are in generic position, and then the rule for calculating the limit moments of ${\bf C}$ are given by Voiculescu's concept of ``freeness.''

\begin{theorem}[Voiculescu \cite{voiculescu91a}]
Let ${\bf A}$ and ${\bf B}$ be two random matrix ensembles of $N\times N$-random matrices,
${\bf A}=({\bf A}_N)_{N\in\NN}$ and ${\bf B}=({\bf B}_N)_{N\in\NN}$, each of them with a limit eigenvalue distribution.
Assume that ${\bf A}$ and ${\bf B}$ are independent (i.e., for each $N\in\NN$, all entries of ${\bf A}_N$ are
independent from all entries of ${\bf B}_N$), and that at least one of them is unitarily invariant
(i.e., for each $N$, the joint distribution of the entries does not change if we conjugate the random matrix
with an arbitrary unitary $N\times N$ matrix).
Then ${\bf A}$ and ${\bf B}$ are asymptotically free in the sense of the following definition.
\end{theorem}

\begin{definition}[Voiculescu \cite{voiculescu85a}]\label{def:freeness-first}
Two random matrix ensembles ${\bf A}=({\bf A}_N)_{N\in\NN}$ and ${\bf B}=({\bf B}_N)_{N\in\NN}$ with limit eigenvalue
distributions are \emph{asymptotically free} if we have for all $p\geq 1$ and all
$n(1),m(1),\dots,n(p)$, $m(p)\geq 1$ that
\begin{multline*}
\lim_{N\to\infty} E\Bigl[\tr \bigl\{ ({\bf A}_N^{n(1)}-\alpha_{n(1)}^A 1)
\cdot ({\bf B}_N^{m(1)}-\alpha_{m(1)}^B 1)
\cdots
({\bf A}^{n(p)}-\alpha_{n(p)}^A 1)\cdot ({\bf B}^{m(p)}-\alpha_{m(p)}^B 1) \bigr\}\Bigr]=0
\end{multline*}
\end{definition}

In essence, asymptotic freeness is actually a rule which allows to calculate all mixed moments in ${\bf A}$ and ${\bf B}$, \ie, all expressions of the form
$$\lim_{N\to\infty} E[\tr({\bf A}^{n(1)}{\bf B}^{m(1)}{\bf A}^{n(2)}{\bf B}^{m(2)}\cdots {\bf A}^{n(p)}{\bf B}^{m(p)})]$$
out of the limit moments of ${\bf A}$ and the limit moments of ${\bf B}$. In particular, this means that
all limit moments of ${\bf A}+{\bf B}$ (which are sums of mixed moments) exist, thus ${\bf A}+{\bf B}$ has a limit
distribution, and are actually determined in terms of the limit moments of ${\bf A}$ and the limit
moments of ${\bf B}$. For more on free probability, including extensions to the setting where the moments do not exist, we refer the reader to \cite{hiai00a,voiculescu:1992,speicher:book,biane03a}.

We now clarify the connection between the operational law of a subclass of algebraic random matrices and the convolution operations of free probability. This will bring into sharp focus how the polynomial method constitutes a framework for computational free probability theory.

\begin{proposition}\label{prop:free convolution}
Let ${\bf A}_{N} \toinp A$ and ${\bf B}_{N} \toinp B$ be two asymptotically free random matrix sequences as in Definition \ref{sec:moments and freeness}.  Then ${\bf A}_{N}+{\bf B}_{N} \toinp A+B$ and ${\bf A}_{N} \times {\bf B}_{N} \toinp AB$ (where the product is defined whenever ${\bf A}_{N} \times {\bf B}_{N}$ has real eigenvalues for every ${\bf A}_{N}$ and ${\bf B}_{N}$) with the corresponding limit eigenvalue density functions, $f_{A+B}$ and $f_{AB}$ given by
\begin{subequations}
\begin{equation}
f_{A+B} = f_{A} \boxplus f_{B}
\end{equation}
\begin{equation}
f_{AB} = f_{A} \boxtimes f_{B} 
\end{equation}
\end{subequations}
where $\boxplus$ denotes free additive convolution and $\boxtimes$ denotes free multiplicative convolution. These convolution operations can be expressed in terms of the R and S transforms as described in Propositions \ref{th:free add} and \ref{th:free multiply} respectively.
\end{proposition}
\begin{proof}
This result appears for density functions with compact support in \cite{voiculescu86a,voiculescu87a}. It was later strengthened to the case of density functions with unbounded support. See \cite{hiai00a} for additional details and references.
\end{proof}

\vspace{0.15in}
In Theorems \ref{th:free rtransform rule} and \ref{th:free stransform rule} we, in effect, showed that the free convolution of algebraic densities produces an algebraic density.  This stated succinctly next.\\

\theorembox{
\begin{corollary}\label{cor:free add}
Algebraic probability distributions form a semi-group under free additive convolution.
\end{corollary}

\begin{corollary}\label{cor:free multiply}
Algebraic distributions with positive semi-definite support form a semi-group under free multiplicative convolution.
\end{corollary}
}

\begin{table}[t]
\centering
\begin{tabular}{|l||c|c|}
\hline\
                                    &                 & \\
{\large{\bf Free additive convolution}} & {\large$f_{A+B} = f_{A} \bp{} f_{B}$}  & {\large$\lrgs{A+B}=\lrgs{A} \bp{r} \lrgs{B}$} \\
& & \\
{\large{\bf Free multiplicative convolution} }& {\large$f_{A\times B} = f_{A} \bt{} f_{B}$ } & {\large$\lsys{A \times B}=\lsys{A} \bt{s} \lsys{B}$}\\
& & \\
\hline
\end{tabular}
\caption{Implicit representation of the free convolution of two algebraic probability densities.}
\label{tab:resultant summary}
\end{table}\index{resultants!free probability}\index{free probability!convolutions}

This establishes a framework for computational free probability theory by identifying the class of distributions for which the free convolution operations produce a ``computable'' distribution.

\subsection{Implicitly encoding the free convolution computations}\index{free convolution!implicit encoding}

The computational framework established relies on being able to implicitly encode free convolution computations as a resultant computation on appropriate bivariate polynomials as in Table \ref{tab:resultant summary}. This leads to the obvious question: Are there other more \textit{effective} ways to implicitly encode free convolution computations? The answer to this rhetorical question will bring into sharp focus the reason why the bivariate polynomial encoding  at the heart of the polynomial method is indispensable for any symbolic computational implementation of free convolution. First, we answer the analogous question about the most effective encoding for classical convolution computations.

Recall that classical convolution can be expressed in terms of the Laplace transform of the distribution function. In what follows, we assume that the distributions have finite moments\footnote{In the general case, tools from complex analysis can be used to extend the argument.}. Hence the Laplace transform can be written as a formal exponential moment generating function. Classical additive and multiplicative convolution of  two distributions produces a distribution whose exponential moment generating function equals  the series (or Cauchy) product and the coefficient-wise (or Hadamard) product of the individual exponential moment generating functions, respectively. Often, however, the Laplace transform of either or both the individual distributions being convolved cannot be written in closed form. The next best thing to do then is to find an implicit way to encode the Laplace transform and to do the convolution computations  via this representation.

When this point of view is adopted, the task of identifying candidate encodings is reduced to finding the class of representations of the exponential generating function that remains closed under the Cauchy and  Hadamard product. Clearly, rational generating functions (see Definition \ref{def:rational gfun}) satisfy this requirement. It is shown in Theorem 6.4.12 \cite[pp.194]{stanley99a}, that D-finite generating functions (see Definition \ref{def:dfinite}) satisfy this requirement as well.

 Proposition \ref{th:dfinite theorem} establishes that all algebraic generating functions (see Definition \ref{def:algebraic gfun}) and by extension, rational generating functions, are also D-finite. However, not all D-finite generating functions are algebraic (see Exercise 6.1 \cite[pp. 217]{stanley99a} for a counter-example) so that algebraic generating functions do not satisfy the closure requirement. Furthermore, from Proposition 6.4.3 and Theorem 6.4.12 in \cite{stanley99a}, if the \textit{ordinary} generating function is D-finite then so is the \textit{exponential} generating function and vice versa. Thus D-finite generating functions are the largest class of generating functions for which classical convolution computations can be performed via an implicit representation. 

In the context of developing a computational framework based on the chosen implicit representation, it is important to consider computability and algorithmic efficiency issues. The class of D-finite functions is well-suited in that regard as well \cite{salvy94a} so that we regard it as the most \textit{effective} class of representations in which the classical convolution computations may be performed implicitly.

However, this class is inadequate for performing free convolution computations implicitly. This is a consequence of the prominent role occupied in this theory  by ordinary generating functions. Specifically, the ordinary formal R and S power series, are obtained from the ordinary moment generating function by functional inversion (or reversion), and are the key ingredients of free additive and multiplicative convolution (see Propositions \ref{prop:free convolution}, \ref{th:free add} and \ref{th:free multiply}).  The task of identifying candidate encodings is thus reduced to finding the class of representations of the ordinary moment generating function that remains closed under addition, the Cauchy product, \textit{and} reversion. D-finite functions only satisfy the first two conditions and are hence unsuitable representations. 

Algebraic functions do, however, satisfy all three conditions. The algorithmic efficiency of computing the resultant (see Section \ref{sec:resultant connection}) justifies our labelling of the bivariate polynomial encoding as the most \textit{effective} way of implicitly encoding free convolution computations. The candidacy of constructibly D-finite generating functions \cite{bergeron90a}, which do not contain the class of D-finite functions but do contain the class of algebraic functions,  merits further investigation since they are closed under reversion, addition and multiplication.  Identifying classes of representations of generating functions for which  \textit{both} the classical and free convolution computations can be performed implicitly and effectively remains an important open problem. 

%\chapter{The polynomial method: Applications}\label{chap:pol method examples}
\section{Applications}\label{sec:pol method examples}
We illustrate the use of the computational techniques developed in Section \ref{sec:pol method computational aspects} with some examples. Documented MATLAB implementation of the polynomial method is available via the RMTool package \cite{raj:rmtool} from \texttt{http://www.mit.edu/\~\/raj/rmtool/}; the examples considered in this article, along with many more, appear there and in \cite{raj:thesis}.

\subsection{The Jacobi random matrix}\index{MANOVA}\index{Jacobi random matrix}\index{Multivariate F matrix}
The Jacobi matrix ensemble is defined in terms of two independent Wishart matrices ${\bf W}_{1}(c_{1})$ and ${\bf W}_{2}(c_{2})$ as ${\bf J}=({\bf I}+{\bf W}_{2}(c_{2})\, {\bf W}_{1}^{-1}(c_{1}))^{-1}$. The subscripts are not to be confused for the size of the matrices. Listing the computational steps needed to generate a realization of this ensemble, as in Table \ref{tab:jacobi example}, is the easiest way to identify the sequence of random matrix operations needed to obtain $\lmzs{J}$.
\begin{table}[h]
\centering
\begin{footnotesize}
\begin{tabular}{|c||l|l|}\hline
Transformation & Numerical MATLAB code & Symbolic MATLAB code\\ \hline
\begin{tabular}{c}
Initialization
\end{tabular} &
\begin{tabular}{@{}l}
 \% Pick n, c1, c2\\
N1=n/c1; N2=n/c2;
\end{tabular} &
\begin{tabular}{@{}l}
\% Define symbolic variables\\
syms m c z;
\end{tabular}\\
&&\\[-1ex]
${\bf A}_1 = {\bf I}$ &
A1 = eye(n,n); &
Lmz1 = m*(1-z)-1;\\
&&\\[-1ex]
${\bf A}_2= {\bf W}_1(c_1) \times {\bf A}_1$ &
\begin{tabular}{@{}l}
G1 = randn(n,N1)/sqrt(N1);\\
W1 = G1*G1';\\
A2 = W1*A1;
\end{tabular} &
Lmz2 = AtimesWish(Lmz1,c1);\\
&&\\[-1ex]
${\bf A}_3 = {\bf A}_2^{-1}$ &
 A3 = inv(A2); &
Lmz3 = invA(Lmz2);\\
&&\\[-1ex]
${\bf A}_4 = {\bf W}_2(c_2) \times {\bf A}_3$ &
\begin{tabular}{@{}l}
G2 = randn(n,N2)/sqrt(N2);\\
W2 = G2*G2';\\
A4 = W2*A3;
\end{tabular} &
Lmz4 = AtimesWish(Lmz3,c2);\\
&&\\[-1ex]
${\bf A}_5 = {\bf A}_4+{\bf I}$ &
A5 = A4+I; &
Lmz5 = shiftA(Lmz4,1);\\
&&\\[-1ex]
${\bf A}_6 = {\bf A}_5^{-1}$ &
A6 = inv(A5); &
Lmz6 = invA(Lmz5);\\ \hline
\end{tabular}
\end{footnotesize}
\caption{Sequence of MATLAB commands for sampling the Jacobi ensemble. The functions used to generate the corresponding bivariate polynomials symbolically are listed in Table \ref{tab:summary bivariate transforms}}.
\label{tab:jacobi example}
\end{table}
\newpage
\noindent We first start off with ${\bf A}_{1}={\bf I}$. The bivariate polynomial that encodes the Stieltjes transform of its eigenvalue distribution function is given by 
\begin{equation}\label{eq:ex1 eq1}
\smallbox{\lmz{1}= (1-z)m-1}.
\end{equation}
For ${\bf A}_{2}={\bf W}_{1}(c_{1})\times {\bf A}_{1}$, we can use (\ref{eq:product pol transform}) to obtain the bivariate polynomial
\begin{equation}\label{eq:ex1 eq2}
\smallbox{\lmz{2} = z\,c_{1}\,{m}^{2}- \left( -c_{1}-z+1 \right) m +1}.
\end{equation}
\noindent For ${\bf A}_{3}={\bf A}_{2}^{-1}$, from (\ref{eq:mobius pol transform}), we obtain the bivariate polynomial
\begin{equation}\label{eq:ex1 eq3}
\smallbox{\lmz{3} = {z}^{2} c_{1}{m}^{2}+ \left( c_{1}\,z+z-1 \right) m+1} .
\end{equation}
\noindent For $A_{4}= {\bf W}_{2}(c_{2}) \times {\bf A}_3$. We can use (\ref{eq:product pol transform}) to obtain the bivariate polynomial
\begin{equation}\label{eq:ex1 eq4}
\smallbox{\lmz{4} = \left(  c_{1}\,{z}^{2}+c_{2}\,z \right) {m}^{2}+ \left( c_{1}\,z+z-1+ c_{2} \right) m +1} .
\end{equation}
\noindent For ${\bf A}_{5}={\bf A}_{4}+{\bf I}$, from (\ref{eq:mobius pol transform}), we obtain the bivariate polynomial
\begin{equation}\label{eq:ex1 eq5}
\smallbox{\lmz{5} =  \left(  \left( z-1 \right) ^{2}c_{1}+c_{2}\, \left( z-1 \right)  \right) {m}^{2}+\left(c_{1}\, \left( z-1 \right) +z-2+ c_{2} \right) m +1}.
\end{equation}
\noindent Finally, for ${\bf J}={\bf A}_{6}={\bf A}_{5}^{-1}$, from (\ref{eq:mobius pol transform}), we  obtain the required bivariate polynomial
\begin{multline}\label{eq:ex1 eq6}
\smallbox{\lmz{J} \equiv \lmz{6} = \left( c_{1}\,z+{z}^{3} c_{1}-2\, c_{1}\,{z}^{2}- c_{2}\,{z}^{3}+ c_{2}\,{z}^{2} \right) {m}^{2}}\\
\smallbox{+ \left( -1+2\,z+ c_{1}-3\, c_{1}\,z+2\, c_{1}\,{z}^{2}+c_{2}\,z-2\,c_{2}\,{z}^{2} \right) m- c_{2}\,z-c_{1}+2+ c_{1}\,z} .
\end{multline}
\noindent Using matrix theoretic arguments, it is clear that the random matrix ensembles ${\bf A}_{3}, \ldots {\bf A}_6$ are defined only when $c_{1}<1$. There will be an atomic mass of weight $(1-1/c_{2})$ at $1$ whenever $c_{2}>1$. The non-atomic component of the distribution will have a region of support $\rocnon{} = (a_{-},a_{+})$. The limiting density function for each of these ensembles can be expressed as
\begin{empheq}[
box=\setlength{\fboxrule}{1pt}\fbox]{equation}\label{eq:ex1 edf}
f_{A_{i}}(x)=\dfrac{\sqrt{(x-a_{-})(a_{+}-x)}}{2\,\pi\,l_{2}(x)} \qquad \textrm{ for } a_{-} < x < a_{+},
\end{empheq}
for $i=2,\ldots,6$, where $a_{-}$, $a_{+}$ , where the polynomials $l_{2}(x)$ are listed in Table \ref{tab:jacobi edf}.
\begin{table}[t]
\centering \footnotesize{
\begin{tabular}{|c|c|c|}
\hline
 & $l_{2}(x)$ & $a_{\pm}$ \\[0.5ex]
\hline \hline
$A_2$ & $x\,c_{1}$ &$(1\pm \sqrt{c_{1}})^{2}$ \\[0.5ex]
\hline
$A_3$ & $x^{2}\,c_{1}$ & $\dfrac{1}{(1\mp \sqrt{c_{1}})^{2}}$\\[1ex]
\hline
$A_4$ & $c_{1}x^{2}+c_{2}x$ & $\dfrac{1+c_{1}+c_{2}-c_{1}c_{2}\pm 2\sqrt{c_{1}+c_{2}-c_{1}c_{2}}}{(1-c_{1})^{2}}$ \\[1ex]
\hline
$A_5$ & $c_{1}(x-1)^{2}+c_{2}(x-1)$ & $\dfrac{c_{1}^{2}-c1+2+c_{2}-c_{1}c_{2}\pm 2\sqrt{c_{1}+c_{2}-c_{1}c_{2}}}{(1-c_{1})^{2}}$\\[1ex]
\hline
$A_6$ & $ \left( c_{1}\,x+{x}^{3} c_{1}-2\, c_{1}\,{x}^{2}- c_{2}\,{x}^{3}+ c_{2}\,{x}^{2} \right)$ & $\dfrac{(1-c_{1})^{2}}{c_{1}^{2}-c1+2+c_{2}-c_{1}c_{2}\mp 2\sqrt{c_{1}+c_{2}-c_{1}c_{2}}}$\\[1ex]
\hline
\end{tabular}
} 
\caption{Parameters for determining the limiting eigenvalue density function using (\ref{eq:ex1 edf}).} 
\label{tab:jacobi edf}
\end{table}
\begin{figure}[t]
 \centering
       \includegraphics[scale=0.80]{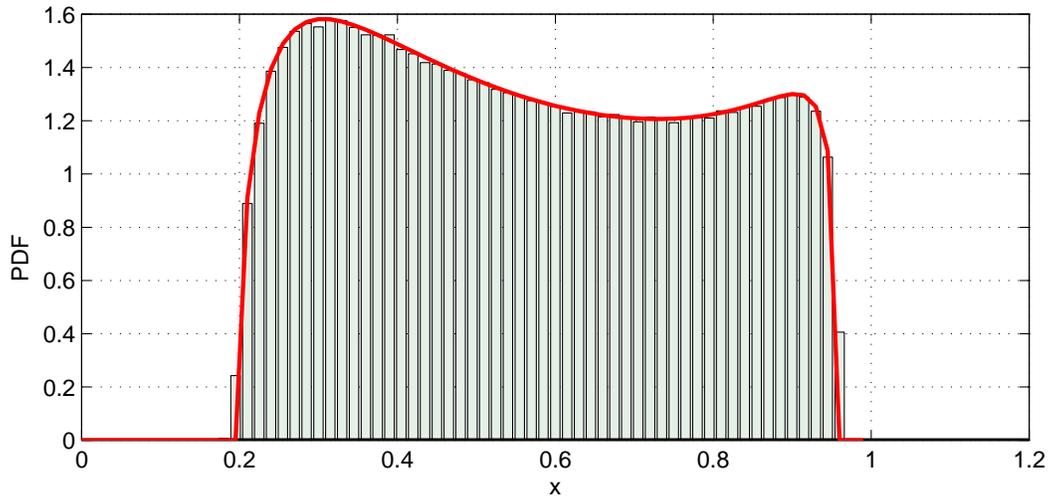}
\label{tab:jacobi random matrix}
\caption{The limiting density (solid line), $f_{A_{6}}(x)$, given by (\ref{eq:ex1 edf}) with $c_{1}=0.1$ and $c_{2}=0.625$ is compared with the normalized histogram of the eigenvalues of a Jacobi matrix generated using the code in Table \ref{tab:jacobi example} over $4000$ Monte-Carlo trials with $n=100$, $N_{1}=n/c_{1} =1000$ and $N_{2}=n/c_{2}=160$.}
\end{figure}
The moments for the general case when $c_{1} \neq c_{2}$ can be enumerated using the techniques described; they will be quite messy. Instead, consider the special case when $c_1=c_2=c$. Using the tools described, the first four terms of the moment series, $\mu(z)=\mu_{J}(z)$, can be computed directly from $\lmuzs{J}$ as
\begin{multline*}
\smallbox{\mu(z)={\frac {1}{2}}+ \left( \frac{1}{8}\,c+{\frac {1}{4}} \right) z+ \left( \frac{3}{16}
\,c+{\frac {1}{8}} \right) {z}^{2}+ \left( \frac{1}{32}\,{c}^{2}+\frac{3}{16}\,c-{
\frac {1}{128}}\,{c}^{3}+{\frac {1}{16}} \right) {z}^{3}}\\ \smallbox{+ \left( -{
\frac {5}{256}}\,{c}^{3}+{\frac {5}{64}}\,{c}^{2}+{\frac {5}{32}}\,c+{
\frac {1}{32}} \right) {z}^{4}+O \left( {z}^{5} \right)}. 
\end{multline*}
The moment generating function satisfies the differential equation
\begin{multline*}
\smallbox{-3\,z+2+zc+ \left( -6\,{z}^{2}+{z}^{3}+10\,z+{z}^{3}{c}^{2}-2\,{z}^{3}
c-4 \right) \mu \left( z \right)} \\
\smallbox{+\left( {z}^{4}-5\,{z}^{3}-2\,{z}^{4}c+8\,{z}^{2}+{z}^{4}{c}^{2}+2\,{z}^{3}c-4\,z-{z}^{3}{c}^{2}
 \right) {\frac {d}{dz}}\mu \left( z \right) =0},
\end{multline*}
with the initial condition $\mu(0) = 1$. The moments $a(n)=M_{n}$ themselves are P-recursive and obtained by the recursion
\begin{multline*}
\smallbox{\left( -2\,c+{c}^{2}+1+ \left( -2\,c+{
c}^{2}+1 \right) n \right) a \left( n \right) + \left(  \left( -5+2\,c
-{c}^{2} \right) n-11+2\,c-{c}^{2} \right) a \left( n+1 \right)} \\
\smallbox{+ \left( 26+8\,n \right) a \left( n+2 \right) + \left( -16-4\,n \right) a \left( n+3 \right)=0}, 
\end{multline*}
with the initial conditions $a(0)=1/2$, $a(1)=1/8\,c+1/4$, and  $a(2)=3/16\,c+1/8$. We can similarly compute the recursion for the free cumulants, $a(n)=K_{n+1},$ as
\[
\smallbox{n{c}^{2}a \left( n \right) +
 \left( 12+4\,n \right) a \left( n+2 \right) =0},
 \]
with the initial conditions $a(0)=1/2$, and $a(1)=1/8\,c$.
\subsection{Random compression of a matrix}
\begin{figure}[t]
\centering
     \includegraphics[scale=0.80]{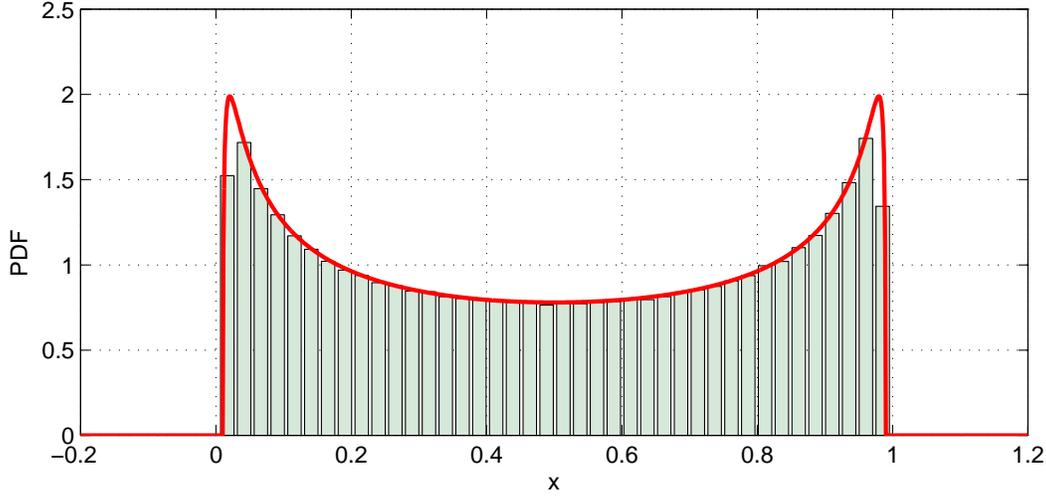}
    \label{fig:compressatom example}
    \caption{The limiting eigenvalue density function (solid line) of the top $0.4N \times 0.4N$ block of a randomly rotated matrix is compared with the experimental histogram collected over $4000$ trials with $N = 200$. Half of the eigenvalues of the original matrix were equal to one while the remainder were equal to zero.}
\end{figure}

\index{random matrices!compression}
\theorembox{
\begin{theorem}\label{th:random compression}
Let ${\bf A}_{N} \mapsto A \in \Palg$. Let ${\bf Q}_{N}$ be an $N \times N$ Haar unitary/orthogonal random matrix independent of ${\bf A}_{N}$. Let ${\bf B}_{n}$ be the upper $n \times n$ block of  ${\bf Q}_{N} {\bf A}_{N} {\bf Q}_{N}^{'}$. Then $${\bf B}_{n} \mapsto B \in \Palg$$ as $n/N \to c$ for $n,N \to \infty$.
\end{theorem}}
\begin{proof}
Let ${\bf P}_{N}$ be an $N \times N$ projection matrix
\[
{\bf P}_{N} \equiv {\bf Q}_{N}
\begin{bmatrix}
{\bf I}_{n} &  \\
 & {\bf 0}_{N-n} \\
\end{bmatrix}
{\bf Q}_{N}^{'}.
\]
By definition, ${\bf P}_{N}$ is an atomic matrix so that ${\bf P}_{N} \to P \in \CM$ as $n/N \to c$ for $n,N \to \infty$. Let $\widetilde{{\bf B}}_{N} = {\bf P}_{N} \times {\bf A}_{N}$. By Corollary \ref{th:free convolution invariant}, $\widetilde{{\bf B}}_{N} \to \widetilde{B} \in \CM$. Finally, from Theorem \ref{th:subspace projection}, we have that ${\bf B}_{n} \to B \in \CM$.  
\end{proof}\\
The proof above provides a recipe for computing the bivariate polynomial $\lmzs{B}$ explicitly as a function of $\lmzs{A}$ and the compression factor $c$. For this particular application, however, one can use first principles \cite{speicher03a} to directly obtain the relationship 
\[
r_{B}(g) = r_{A}(c\,g),
\]
expressed in terms of the R transform. This translates into the operational law
\begin{empheq}[
box=\setlength{\fboxrule}{1pt}\fbox]
{equation}\label{eq:compression shortcut}
\lrg{B} = \lrgs{A}(r,c\,g). 
\end{empheq}
\noindent \textbf{Example}: Consider the atomic matrix ${\bf A}_{N}$ half of whose eigenvalues are equal to one while the remainder are equal to zero. Its eigenvalue distribution function is given by (\ref{eq:edf example}). From the bivariate polynomial, $\lrgs{A}$ in Table \ref{tab:example 1 summary} and (\ref{eq:compression shortcut}) it can be show that the limiting eigenvalue distribution function of ${\bf B}_{n}$, constructed from ${\bf A}_{N}$ as in Theorem \ref{th:random compression}, is encoded by the polynomial
\[
\smallbox{\lmzs{B}  = \left( -2\,c{z}^{2}+2\,cz \right) {m}^{2}- \left( -2\,c+4\,cz+1-2\,z \right) m-2\,c+2},
\]
where $c$ is the limiting compression factor. Poles occur at $z=0$ and $z=1$. The leading terms of the Puiseux expansion of the two branches about the poles at $z=z_{0}$ are
\[
\smallbox{\left\{  \left( {\frac {z-z_{0}}{-2\,c+4\,{c}^{2}}}+{\frac {1-2\,c}{2c}
} \right)\frac{1}{z-z_{0}},{\frac {2\,c-2}{-1+2\,c}} \right\}}.
\]
It can be easily seen that when $c>1/2$, the Puiseux expansion about the poles $z=z_{0}$ will correspond to an atom of weight $w_{0}=(2c-1)/2c$. Thus the limiting eigenvalue distribution function has density
\begin{empheq}[
box=\setlength{\fboxrule}{1pt}\fbox]{equation}\label{eq:compression answer}
\smallbox{f_{B}(x) = \max \left(\frac{2c-1}{2c},0\right) \delta(x) + \frac{1}{\pi}\frac{\sqrt{(x-a_{-})(a_{+}-x)}}{2\,xc-2\,c{x}^{2}} I_{[a_{-},a_{+}]}+ \max \left(\frac{2c-1}{2c},0\right) \delta(x-1)},
\end{empheq}
where $a_{\pm} = 1/2\pm \sqrt {-{c}^{2}+c}$. Figure \ref{fig:compressatom example} compares the theoretical prediction in (\ref{eq:compression answer}) with a Monte-Carlo experiment for $c=0.4$.
From the associated bivariate polynomial
\begin{equation*}
\lmuzs{B}\equiv \left( -2\,c+2\,cz \right) {\mu}^{2}+ \left( z-2-2\,cz+4\,c \right) {\mu}-2\,c+2,
\end{equation*} 
we obtain two series expansions whose branches tend to zero. The first four terms of the series are given by
\begin{equation}\label{eq:true series}
\smallbox{1+{\frac {1}{2}}z+ \frac {1+c}{4} {z}^{2}+ \frac {3+c}{8} {z}^{3}+O \left( {z}^{4} \right)},
\end{equation}
and,
\begin{equation}\label{eq:false series}
\smallbox{{\frac {c-1}{c}}+{\frac {c-1}{2c}}z-{\frac { \left( c-1 \right)  \left( -2+c \right) }{4c}}{z}^{2}-{\frac{ \left( c-1 \right)  \left( 3\,c-4 \right) }{8c}}{z}^{3}+O\left( {z}^{4} \right)}, 
\end{equation}
\noindent respectively. Since $c \leq 1$, the series expansion in (\ref{eq:false series}) can be eliminated since $\mu(0): = \int dF^{B}(x) = 1$. Thus the coefficients of the series in (\ref{eq:true series}) are the correct moments of the limiting eigenvalue distribution. A recursion for the moments can be readily derived using the techniques developed earlier. 
\subsection{Free additive convolution of equilibrium measures}\index{equilibrium measure!convolution}
\begin{figure}[t]
 \centering
     \includegraphics[scale=0.80]{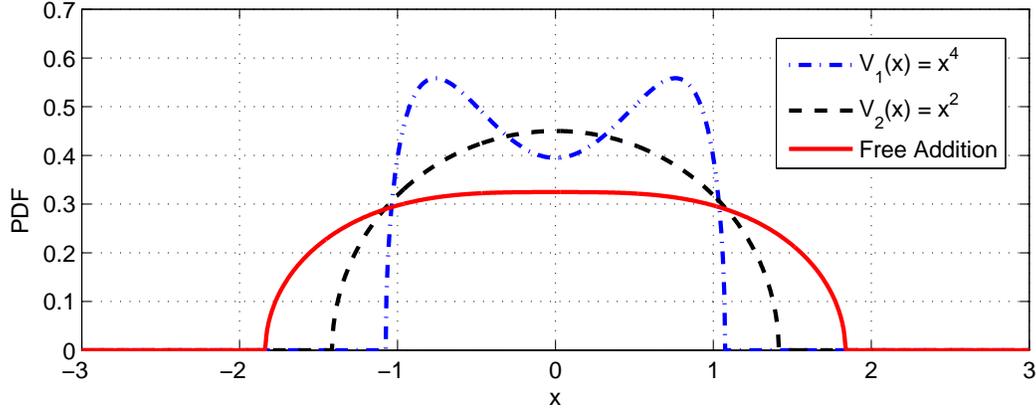}
       \label{fig:eqmeasure example}
\caption{Additive convolution of equilibrium measures corresponding to potentials $V_{1}(x)$ and $V_{2}(x)$.}   
\end{figure}
Equilibrium measures are a fascinating topic within random matrix theory. They arise in the context of research that examines why very general random models for random matrices exhibit universal behavior in the large matrix limit. Suppose we are given a potential $V(x)$ then we consider a sequence of Hermitian, unitarily invariant random matrices ${\bf A}_{N}$, the joint distribution of whose elements is of the form
\[
P({\bf A}_{N}) \propto \exp\left(-N \,{\rm Tr}\, V({\bf A}_{N})\right) d{\bf A}_{N},
\]
where $d{\bf A}_{N} = \prod_{i \leq j} (d{\bf A}_{N})_{ij}$. The equilibrium measure, when it exists, is the unique probability distribution function that minimizes the logarithmic energy (see \cite{deift_book} for additional details). The resulting equilibrium measure depends explicitly on the potential $V(x)$ and can be explicitly computed for some potentials. In particular, for potentials of the form $V(x) = t\,x^{2m}$, the Stieltjes transform of the resulting equilibrium measure is an algebraic function \cite[Chp. 6.7, pp. 174-175]{deift_book} so that the equilibrium measure is an algebraic distribution. Hence we can formally investigate the additive convolution of equilibrium measures corresponding to two different potentials. For $V_{1}(x) = x^{2}$, the equilibrium measure is the (scaled) semi-circle distribution encoded by the bivariate polynomial
\[
\smallbox{\lmzs{A} \equiv m^2+2\,m\,z+2}.
\]
For $V_{2}(x)=x^{4}$, the equilibrium measure is encoded by the bivariate polynomial
\[
\smallbox{\lmzs{B} \equiv 1/4\,{m}^{2}+m{z}^{3}+{z}^{2}+2/9\,\sqrt {3}}.
\]
Since ${\bf A}_{N}$ and ${\bf B}_{N}$ are unitarily invariant random matrices, if ${\bf A}_{N}$ and ${\bf B}_{N}$ are independent, then the limiting eigenvalue distribution function of ${\bf C}_{N} = {\bf A}_{N} + {\bf B}_{N}$ can be computed from $\lmzs{A}$ and $\lmzs{B}$. The limiting eigenvalue density function $f_{C}(x)$ is the free additive convolution of $f_{A}$ and $f_{B}$. The \matlab command \texttt{LmzC = AplusB(LmzA,LmzB);} will produce the bivariate polynomial
\[
\smallbox{\lmzs{C} = -9\,{m}^{4}-54\,{m}^{3}z+ \left( -108\,{z}^{2}-36 \right) {m}^{2}- \left( 72\,{z}^{3}+72\,z \right) m-72\,{z}^{2}-16\,\sqrt {3}}.
\]
Figure \ref{fig:eqmeasure example} plots the probability density function for the equilibrium measure for the potentials $V_{1}(x) = x^2$ and $V_{2}(x) = x^{4}$ as well as the free additive convolution of these measures. The interpretation  of the resulting measuring in the context of potential theory is not clear. The matrix ${\bf C}_{N}$ will no longer be unitarily invariant so it is pointless to look for a potential $V_{3}(x)$ for which $F^{C}$ is an equilibrium measure. The tools and techniques developed in this article might prove useful in further explorations.

\subsection{Algebraic sample covariance matrices}
The (broader)  class of algebraic Wishart sample covariance matrices for which this framework applies is described next.

\theorembox{
\begin{theorem}\label{th:spatio temporal wishart}
  Let ${\bf A}_{n} \toinp A \in \CM$,  and ${\bf B}_{N} \toinp B \in \CM$ be algebraic covariance matrices with ${\bf G}_{n,N}$ denoting an $n \times N$ (pure) Gaussian random matrix (see Definition \ref{def:gaussianlike matrix}).  Let ${\bf X}_{n,N} = {\bf A}_{n}^{1/2} {\bf G}_{n,N} {\bf B}_{N}^{1/2}$. Then $${\bf S}_{n} = {\bf X}_{n,N} {\bf X}_{n,N}^{'} \toinp S \in \CM,$$
as $n,N \to \infty$ and $c_{N} = n/N \to c$.
\end{theorem}} 
\begin{proof}
Let ${\bf Y}_{n,N} \equiv {\bf G}_{n,N} {\bf B}_{N}^{1/2}$, ${\bf T}_{n} \equiv {\bf Y}_{n,N} {\bf Y}_{n,N}'$ and $\widetilde{{\bf T}}_{N} = {\bf Y}_{n,N}' {\bf Y}_{n,N}$. Thus ${\bf S}_{n} = {\bf A}_{n} \times {\bf T}_{n} \equiv {\bf A}_{n}^{1/2} {\bf T}_{n} {\bf A}_{n}^{1/2}$. The matrix ${\bf T}_{n}$, as defined, is invariant under orthogonal/unitary transformations, though the matrix $\widetilde{{\bf T}}_{N}$ is not. Hence, by Corollary \ref{th:free convolution invariant}, and since ${\bf A}_{n} \mapsto A \in \CM$, ${\bf S}_{n} \mapsto S \in \CM$ whenever ${\bf T}_{n} \mapsto T \in \CM$. 

From Theorem \ref{th:transpose}, ${\bf T}_{n} \mapsto T \in \CM$ if $\widetilde{{\bf T}}_{N} \mapsto \widetilde{T} \in \CM$. The matrix $\widetilde{{\bf T}}_{N} = {\bf B}_{N}^{1/2} {\bf G}_{n,N}'{\bf G}_{n,N} {\bf B}_{N}^{1/2}$ is clearly algebraic by application of Corollary \ref{th:free convolution invariant} and Theorem \ref{th:simple deterministic} since ${\bf B}_{N}$ is algebraic and ${\bf G}_{n,N}'{\bf G}_{n,N}$ is algebraic and unitarily invariant. 
\end{proof}\\

The proof of Theorem \ref{th:spatio temporal wishart} provides us with a recipe for computing the polynomials that encode the limiting eigenvalue distribution of ${\bf S}$ in the far more general situation where the observation vectors are modelled as samples of a multivariate Gaussian with spatio-temporal correlations. The limiting eigenvalue distribution of ${\bf S}$ depends on the limiting (algebraic) eigenvalue distributions of ${\bf A}$ and ${\bf B}$ and may be called using the \texttt{AtimesWishtimesB} function in the RMTool \cite{raj:rmtool} package. See \cite{raj07a} for the relevant code.

\subsection{Other applications}
There is often a connection between well-known combinatorial numbers and random matrices. For example, the even moments of the Wigner matrix are the famous Catalan numbers. Similarly, if ${\bf W}_{N}(c)$ denotes the Wishart matrix with parameter $c$, other combinatorial correspondences can be easily established using the techniques developed. For instance, the limiting moments of ${\bf W}_{N}(1)-{\bf I}_{N}$ are the Riordan numbers, the large Schr\"{o}der numbers correspond to the limiting moments of $2{\bf W}_{N}(0.5)$ while the small Schr\"{o}der numbers are the limiting moments of $4{\bf W}_{N}(0.125)$. Combinatorial identities along the lines of those developed in \cite{dumitriu03a} might result from these correspondences. \\ 
\section{Some open problems}\label{sec:open problems}
\begin{itemize}
\item{If we are given a single realization of an $N \times N$ sized algebraic random matrix ${\bf A}_{N}$, is it possible to reliably infer the minimal bivariate polynomial $\lmzs{}$, \ie, with combined degree $\Dm+\Dz$ as small as possible, that encodes its limiting eigenvalue distribution?}
\item{This is closely related to the following problem. Define the set of ``admissible'' real-valued coefficients $c_{ij}$ for $0 \leq i \leq \Dm$ and $0 \leq j \leq \Dz$. Here admissibility implies that (a branch of) a solution $m(z)$ of the equation $\lmz{A}:=\sum_{i=0}^{\Dm} \sum_{j=0}^{\Dz} c_{ij} u^{i} v^{j}=0$ is globally the Stieltjes transform  of a positive probability distribution.} 
\end{itemize}
\section*{Acknowledgements}
Our work has greatly benefitted from interactions with colleagues. It is our pleasure to acknowledge their contribution. Jack Silverstein kindled the first author's involvement in the subject with a manuscript of unpublished work and several email discussions.   Richard Stanley introduced us to D-finite series while we learned about resultants from Pablo Parillo. We particularly wish to thank the anonymous referees for carefully proof-reading the original manuscript and challenging us to improve the organization of the paper. 

We owe Roland Speicher thanks for numerous conversations, feedback, encouragement and for his patience in answering our queries. Our enthusiasm for free probability can be directly traced to his inspiring exposition in his lectures we have attended and the survey papers posted on his website. So much so, that when we sat down to write Section \ref{sec:moments and freeness}, we decided that it was pointless attempting to emulate the clarity of his writing. We are grateful to Roland for allowing us to borrow his writing in Section \ref{sec:moments and freeness}.

The authors were supported by a National Science Foundation grant DMS-0411962 and an Office of Naval Research Special Post-Doctoral Award N00014-07-1-0269. 

\bibliographystyle{siam}
\bibliography{randbib}

\def\cprime{$'$} \def\cprime{$'$}
\begin{thebibliography}{10}

\bibitem{akhiezer65a}
{\sc N.~I. Akhiezer}, {\em The classical moment problem and some related
  questions in analysis}, Hafner Publishing Co., New York, New York, 1965.
\newblock Translated by N. Kemmer.

\bibitem{akritas93a}
{\sc A.~G. Akritas}, {\em Sylvester's forgotten form of the resultant},
  Fibonacci Quart., 31 (1993), pp.~325--332.

\bibitem{anderson06c}
{\sc G.~Anderson and O.~Zeitouni}, {\em A law of large numbers for finite-range
  dependent random matrices}.
\newblock \texttt{http://arxiv.org/abs/math/0609364}, September 2006.
\newblock Preprint.

\bibitem{bai95a}
{\sc Z.~D. Bai and J.~W. Silverstein}, {\em On the empirical distribution of
  eigenvalues of a class of large dimensional random matices}, Journal of
  Multi. Analysis, 54 (1995), pp.~175--192.

\bibitem{bergeron90a}
{\sc F.~Bergeron and C.~Reutenauer}, {\em Combinatorial resolution of systems
  of differential equations. {III}. {A} special class of differentially
  algebraic series}, European J. Combin., 11 (1990), pp.~501--512.

\bibitem{biane03a}
{\sc P.~Biane}, {\em Free probability for probabilists}, in Quantum probability
  communications, Vol. XI (Grenoble, 1998), QP-PQ, XI, World Sci. Publishing,
  River Edge, NJ, 2003, pp.~55--71.

\bibitem{billingsley99a}
{\sc P.~Billingsley}, {\em Convergence of probability measures}, Wiley Series
  in Probability and Statistics: Probability and Statistics, John Wiley \& Sons
  Inc., New York, second~ed., 1999.
\newblock A Wiley-Interscience Publication.

\bibitem{collins05a}
{\sc B.~Collins}, {\em Product of random projections, {J}acobi ensembles and
  universality problems arising from free probability}, Probab. Theory Related
  Fields, 133 (2005), pp.~315--344.

\bibitem{deift98a}
{\sc P.~Deift, T.~Kriecherbauer, and K.~T.-R. McLaughlin}, {\em New results on
  the equilibrium measure for logarithmic potentials in the presence of an
  external field}, J. Approx. Theory, 95 (1998), pp.~388--475.

\bibitem{deift99a}
{\sc P.~Deift, T.~Kriecherbauer, K.~T.-R. McLaughlin, S.~Venakides, and
  X.~Zhou}, {\em Uniform asymptotics for polynomials orthogonal with respect to
  varying exponential weights and applications to universality questions in
  random matrix theory}, Comm. Pure Appl. Math., 52 (1999), pp.~1335--1425.

\bibitem{deift_book}
{\sc P.~A. Deift}, {\em Orthogonal polynomials and random matrices: a
  {R}iemann-{H}ilbert approach}, vol.~3 of Courant Lecture Notes in
  Mathematics, New York University Courant Institute of Mathematical Sciences,
  New York, 1999.

\bibitem{silverstein04a}
{\sc W.~B. Dozier and J.~W. Silverstein}, {\em On the empirical distribution of
  eigenvalues of large dimensional information-plus-noise type matrices}.
\newblock {\tt http://www4.ncsu.edu/\~\/jack/infnoise.pdf}, 2004.

\bibitem{dumitriu03a}
{\sc I.~Dumitriu and E.~Rassart}, {\em Path counting and random matrix theory},
  Electronic Journal of Combinatorics, 7 (2003).
\newblock R-43.

\bibitem{flajolet01a}
{\sc P.~Flajolet and R.~Sedgewick}, {\em Analytic combinatorics: Functional
  equations, rational and algebraic functions}, Research Report 4103, INRIA,
  2001.
\newblock {\tt http://algo.inria.fr/flajolet/Publications/FlSe01.pdf}.

\bibitem{hiai00a}
{\sc F.~Hiai and D.~Petz}, {\em The semicircle law, free random variables and
  entropy}, vol.~77, {American Mathematical Society}, 2000.

\bibitem{horn91a}
{\sc R.~A. Horn and C.~R. Johnson}, {\em Topics in matrix analysis}, Cambridge
  University Press, Cambridge, 1991.

\bibitem{kuijlaars00a}
{\sc A.~Kuijlaars and K.~T.-R. McLaughlin}, {\em Generic behavior of the
  density of states in random matrix theory and equilibrium problems in the
  presence of real analytic external fields}, Commun. Pure Appl. Math., 53
  (2000), pp.~736--785.

\bibitem{marcenko67a}
{\sc V.~A. Mar{\v{c}}enko and L.~A. Pastur}, {\em Distribution of eigenvalues
  in certain sets of random matrices}, Mat. Sb. (N.S.), 72 (114) (1967),
  pp.~507--536.

\bibitem{mckay81a}
{\sc B.~D. McKay}, {\em The expected eigenvalue distribution of a large regular
  graph}, Linear Algebra And Its Applications, 40 (1981), pp.~203--216.

\bibitem{raj:thesis}
{\sc R.~R. Nadakuditi}, {\em {A}pplied {S}tochastic {E}igen-{A}nalysis}, PhD
  thesis, {M}assachusetts {I}nstitute of {T}echnology, February 2007.
\newblock {D}epartment of {E}lectrical {E}ngineering and {C}omputer {S}cience.

\bibitem{speicher:book}
{\sc A.~Nica and R.~Speicher}, {\em {L}ectures on the {C}ombinatorics of {F}ree
  {P}robability}, {L}ondon {M}athematical {S}ociety {L}ecture {N}ote {S}eries,
  {N}ew {Y}ork: {C}ambridge {U}niversity {P}ress, 2006.

\bibitem{raj:rmtool}
{\sc N.~R. Rao}, {\em {RMT}ool: A random matrix and free probability calculator
  in {MATLAB}}.
\newblock \texttt{http://www.mit.edu/\~\/raj/rmtool/}.

\bibitem{raj07a}
{\sc N.~R. Rao and A.~Edelman}, {\em The polynomial method for the eigenvectors
  of random matrices}.
\newblock Preprint.

\bibitem{salvy94a}
{\sc B.~Salvy and P.~Zimmermann}, {\em Gfun: a {M}aple package for the
  manipulation of generating and holonomic functions in one variable}, ACM
  Trans. on Math. Software, 20 (1994), pp.~163--177.

\bibitem{silverstein85a}
{\sc J.~W. Silverstein}, {\em The limiting eigenvalue distribution of a
  multivariate {F} matrix}, SIAM Journal on Math. Anal., 16 (1985),
  pp.~641--646.

\bibitem{silverstein95b}
\leavevmode\vrule height 2pt depth -1.6pt width 23pt, {\em Strong convergence
  of the empirical distribution of eigenvalues of large dimensional random
  matrices}, J. of Multivariate Anal., 55(2) (1995), pp.~331--339.

\bibitem{silverstein95a}
{\sc J.~W. Silverstein and S.-I. Choi}, {\em Analysis of the limiting spectral
  distribution of large-dimensional random matrices}, J. Multivariate Anal., 54
  (1995), pp.~295--309.

\bibitem{speicher97a}
{\sc R.~Speicher}, {\em Free probability theory and non-crossing partitions},
  S\'em. Lothar. Combin., 39 (1997), pp.~Art.\ B39c, 38 pp.\ (electronic).

\bibitem{speicher03a}
\leavevmode\vrule height 2pt depth -1.6pt width 23pt, {\em Free probability
  theory and random matrices}, in Asymptotic combinatorics with applications to
  mathematical physics (St. Petersburg, 2001), vol.~1815 of Lecture Notes in
  Math., Springer, Berlin, 2003, pp.~53--73.

\bibitem{stanley99a}
{\sc R.~P. Stanley}, {\em Enumerative combinatorics. {V}ol. 2}, vol.~62 of
  Cambridge Studies in Advanced Mathematics, Cambridge University Press,
  Cambridge, 1999.
\newblock With a foreword by Gian-Carlo Rota and appendix 1 by Sergey Fomin.

\bibitem{sturmfels97a}
{\sc B.~Sturmfels}, {\em Introduction to resultants}, in Applications of
  computational algebraic geometry (San Diego, CA, 1997), vol.~53 of Proc.
  Sympos. Appl. Math., Amer. Math. Soc., Providence, RI, 1998, pp.~25--39.

\bibitem{tulino04a}
{\sc A.~M. Tulino and S.~Verd{\'u}}, {\em Random matrices and wireless
  communications}, Foundations and Trends in Communications and Information
  Theory, 1 (2004).

\bibitem{voiculescu85a}
{\sc D.~Voiculescu}, {\em Symmetries of some reduced free product {$C\sp
  \ast$}-algebras}, in Operator algebras and their connections with topology
  and ergodic theory (Bu\c steni, 1983), vol.~1132 of Lecture Notes in Math.,
  Springer, Berlin, 1985, pp.~556--588.

\bibitem{voiculescu86a}
\leavevmode\vrule height 2pt depth -1.6pt width 23pt, {\em Addition of certain
  noncommuting random variables}, J. Funct. Anal., 66 (1986), pp.~323--346.

\bibitem{voiculescu87a}
\leavevmode\vrule height 2pt depth -1.6pt width 23pt, {\em Multiplication of
  certain noncommuting random variables}, J. Operator Theory, 18 (1987),
  pp.~223--235.

\bibitem{voiculescu91a}
\leavevmode\vrule height 2pt depth -1.6pt width 23pt, {\em Limit laws for
  random matrices and free products}, Invent. Math., 104 (1991), pp.~201--220.

\bibitem{voiculescu:1992}
{\sc D.~V. Voiculescu, K.~J. Dykema, and A.~Nica}, {\em Free random variables},
  vol.~1 of CRM Monograph Series, American Mathematical Society, Providence,
  RI, 1992.

\bibitem{wigner55a}
{\sc E.~P. Wigner}, {\em Characteristic vectors of bordered matrices with
  infinite dimensions}, Annals of Math., 62 (1955), pp.~548--564.

\end{thebibliography}
\end{document}